\documentclass[
12pt]{article}
\usepackage{epsfig}
\usepackage{amssymb,mathrsfs,amsmath}
\usepackage{color}
\numberwithin{equation}{section}
\newtheorem{theo}{Theorem}[section]
\newtheorem{defi}[theo]{Definition}
\newtheorem{rem}[theo]{\it Remark}
\newtheorem{prop}[theo]{Proposition}
\newtheorem{cor}[theo]{Corollary}
\newtheorem{lemma}[theo]{Lemma}

\def\pp{{\rm I\hspace{-1.5pt}P}}     
\def\ee{\mathbb{E}}
\def\dd{{\rm d}}
\def\dB{\int\hspace{-10pt}\int_{{\cal Q}_+}\hspace{-10pt}B(\dd y,\dd s)}

\def\Co{C_c^\infty\hspace{-3pt}\left(\rule{0pt}{10pt}\right.
        \hspace{-3pt}(x_0,\infty)\hspace{-3pt}\left.\rule{0pt}{10pt}\right)}

\def\Box{\rule{6pt}{6pt}}
\def\LL{\mbox{\boldmath$Lin$}}
\newcommand*{\fatdot}{ \raisebox{-1pt}{\makebox[.5ex]{\textbf{$\cdot$}}} }
\def\NN{\mathbb{N}}
\def\NR{\mathbb{R}}

\newcommand{\ind}{{\bf 1}}
\begin{document}
\title{On the spatial dynamics\\ of the solution to the stochastic heat equation}
\author{Sigurd Assing and James Bichard\\
        Department of Statistics, The University of Warwick\\
        e-mail: s.assing@warwick.ac.uk}
\date{}
\maketitle
\begin{abstract}
We consider the solution of 
$\partial_t u=\partial_x^2 u+\partial_x\partial_t B,\,(x,t)\in\NR\times(0,\infty)$,
subject to $u(x,0)=0,\,x\in\NR$, where $B$ is a Brownian sheet.
We show that $u$ also satisfies
$\partial_x^2 u +[\,(-\partial_t^2)^{1/2}+\sqrt{2}\partial_x(-\partial_t^2)^{1/4}\,]\,u^a=
\partial_x\partial_t{\tilde B}$ in $\NR\times(0,\infty)$
where $u^a$ stands for the extension of $u(x,t)$ to $(x,t)\in\NR^2$ which is antisymmetric in $t$
and $\tilde{B}$ is another Brownian sheet.
The new SPDE allows us to prove the strong Markov property of the pair $(u,\partial_x u)$
when seen as a process indexed by $x\ge x_0$, $x_0$ fixed, taking values in a state space of functions in $t$.
The method of proof is based on enlargement of filtration and we discuss
how our method could be applied to other quasi-linear SPDEs.
\end{abstract}
%
\noindent
{\large KEY WORDS}\hspace{0.4cm}
stochastic partial differential equation, enlargement of filtration,\\
Brownian sheet, Gaussian analysis

\vspace{10pt}
\section{Introduction}
When studying stochastic partial differential equations (SPDEs)
one has to understand the behaviour of multi-parameter random fields
$u({\bf x}),\,{\bf x}\in{\cal Q}$, where ${\cal Q}\subseteq\NR^d$ is a given domain.
A first and of course important question deals with the possibility of a
Markovian `behaviour' of such a field. 
Since L\'{e}vy's sharp Markov property (see \cite{L1945})
already fails in the case of multi-parameter Brownian sheets,
the only comprehensive Markovian `behaviour' one can hope for is the so-called
germ Markov property--the reader is referred to the early papers
\cite{K1961,McK1963} and in particular to \cite{K1970}
for a good introduction to this concept.

There is an early paper by Y.A.\ Rozanov \cite{R1977} on the
Markovian `behaviour' of SPDEs and then there are three influential papers
\cite{D-M1992,D-MN1994,NP1994}
on the germ Markov property of solutions to SPDEs of type
$${\cal L}u\,=\,\eta\,+\,f(u)\quad\mbox{in}\quad{\cal Q}$$
where ${\cal L}$ is a linear partial differential operator
of elliptic or parabolic type and $\eta$ stands for a multi-parameter noise.
The method applied in all three papers consisted in, first,
establishing the germ Markov property for the solution of the linear equation
$${\cal L}u\,=\,\eta\quad\mbox{in}\quad{\cal Q}$$
and, second, getting it for the drift-perturbed equation by a
Kusuoka or Girsanov transform. It should be mentioned that
\cite{AFN1995} provides another useful method for the second step.

The main method for the first step is usually based on the paper \cite{P1971}
which was later improved by H.\ K\"unsch \cite{K1979}.
For example, the more recent paper \cite{BK2008} on the germ Markov property
of the solution of a linear stochastic heat equation is still about
checking the conditions stated in \cite{P1971,K1979} which can be demanding.

However, 
the purpose of our paper is to \underline{refine} this first step in the following sense:
study a more specific Markovian `behaviour' of solutions of linear SPDEs.
Our working example is the stochastic heat equation with additive space-time white noise.
But all calculations are based on only two ingredients:
\begin{itemize}
\item
a Green's function for $\quad{\cal L}u\,=\,\eta\quad\mbox{in}\quad{\cal Q}$;
\item
the covariance of a Gaussian noise $\eta$.
\end{itemize}
So, following our scheme of calculations but using a different Green's function or covariance
would produce similar results with respect to other linear SPDEs.
The explicit calculations are involved and will be different in other cases. 
That's why we have to restrict ourselves to the case of an important example
in order to show how the method works in the very detail.
However, at the end of this introduction, we list the tasks to be dealt with
when treating another SPDE.

We now explain what we mean by a {\it Markovian `behaviour' more specific than the 
germ Markov property}. A random field $u({\bf x}),\,{\bf x}\in{\cal Q}$, satisfies the 
germ Markov property if
$\sigma\{u({\bf y}):{\bf y}\in A\}$ and $\sigma\{u({\bf y}):{\bf y}\in A^c\}$
are conditionally independent given the germ-$\sigma$-algebra
$$germ(\partial A)\,\stackrel{\mbox{\tiny def}}{=}\,
\bigcap_{\{\partial A\subset O:O\,open\,in\,{\cal Q}\}}
\sigma\{u({\bf y}):{\bf y}\in O\}$$
for any Borel set $A\in{\cal Q}$.
This type of Markov property is `directionless' with respect to the 
$d$-dimensional domain ${\cal Q}$. But often it is desirable to emphasise 
a direction in $\NR^d$ and to study the behaviour of an SPDE along this direction.
In the case of parabolic SPDEs, a natural direction to emphasise is 
the direction of `time' we denote by $t$ in what follows. 
Many solutions $u({\bf x},t)$ of parabolic SPDEs are even constructed as Markov processes
$t\mapsto u(\cdot,t)$ taking values in a function space.
Hence, along the direction of time, these solutions satisfy a sharp Markov property
with an associated martingale problem.

But we want to be able to pick other directions with respect to the space-variable
${\bf x}=(x_1,\dots,x_{d-1})$ of $u({\bf x},t)$, for example, the direction of $x_1$.
Assume we would already know that $u({\bf x},t),\,({\bf x},t)\in{\cal Q}$, satisfies
the germ Markov property. Then, the process $x_1\mapsto u(x_1,\cdot)$
is at least Markovian relative to
$$germ\left(\rule{0pt}{10pt}\right.
(\{x_1\}\times\NR^{d-1})\cap{\cal Q}
\left.\rule{0pt}{10pt}\right).$$
But this does \underline{not} give an associated martingale problem
with martingales indexed by $x_1$ hence many useful probabilistic methods
\underline{cannot} be applied.
So one wants to know if there is a $\sigma$-algebra included in
$germ\left(\rule{0pt}{10pt}\right.
(\{x_1\}\times\NR^{d-1})\cap{\cal Q}
\left.\rule{0pt}{10pt}\right)$
so that $x_1\mapsto u(x_1,\cdot)$ is still Markovian relative to this smaller
$\sigma$-algebra but also satisfies an associated martingale problem.
And, because we are dealing with SPDEs, it is very likely that 
a $\sigma$-algebra generated
by certain partial derivatives of $u({\bf x},t)$ with respect to $x_1$
would serve the purpose.

Our main result, Theorem \ref{strMtheo},
states that the above can be achieved in the case of the
stochastic heat equation 
$(\partial^2_x-\partial_t)u=-\partial_x\partial_t B,\,(x,t)\in\NR\times(0,\infty)$,
driven by a Brownian sheet $B$.
We find a martingale problem for the pair of processes
$(u(x,\cdot),\partial_x u(x,\cdot)),\, x\ge x_0$,
which is given by an unbounded operator we explicitly calculate
using the technique of enlargement of filtration.

The need for an enlargement of filtration in this context
can be considered the key idea of this paper.
The problem is explained in Section \ref{results}--the reader is referred to 
the second equation of (\ref{basic SDE}). 
The observation is that, for a given test function $h$,
there is \underline{no} filtration such that
$x\mapsto U(x,h')$ is adapted with respect to this filtration and
$x\mapsto W_{x-x_0}(h)$ is a Wiener process with respect to this filtration.
Enlargement of filtration solves this problem subject to a drift correction.
But the new drift requires test functions which are less regular than
the original test functions $h$.
As a consequence one has to discuss the regularity of all involved processes
very carefully.

We are then able to derive, by showing the uniqueness of the martingale problem,
the strong Markov property of $u(x,\cdot),\,x\ge x_0$,
with respect to the natural filtration generated by 
$(u(x,\cdot),\partial_x u(x,\cdot)),\, x\ge x_0$.

As explained earlier, the same method could be used to find interesting
martingale problems associated with other linear SPDEs or even
drift-perturbed linear SPDEs, by applying Girsanov's transform for example,
the latter being beyond the scope of this paper.

The organisation of the paper is as follows. In Section \ref{notation}
we list notation which is crucial for the understanding of the paper.
Section \ref{results} is a combination of results and further explanations 
which fully describes our method  and can be summarised by:
\begin{itemize}
\item
choose a direction along which one wishes to study the solution $u$ of 
a stochastic partial differential equation
${\cal L}u\,=\,\eta$ in ${\cal Q}$ subject to \underline{given} boundary data;
\item
describe the dynamics of ${\cal L}u\,=\,\eta$ in ${\cal Q}$
along the chosen direction---see (\ref{system});
\item
find the regularity of all involved partial derivatives---see Proposition \ref{stateSpace};
\item
find a correction $\varrho$ of $\eta$ such that ${\cal L}u$
and $\tilde{\eta}=\eta-\varrho$ are both adapted with respect to a filtration
along the chosen direction and that the probability law of $\tilde{\eta}$
is the same as the law of $\eta$---see Proposition \ref{new wipro}, 
Lemma \ref{identify rho}, Remark \ref{details of method}(ii);
\item
describe $\varrho$ as a functional of the 
solution $u$---see Proposition \ref{new drift}, Theorem \ref{new SPDE};
\item
show uniqueness of the martingale problem along the chosen direction
which is associated with the 
new equation
${\cal L}u-\varrho(u)\,=\,\tilde{\eta}$---see Theorem \ref{strMtheo}.
\end{itemize}
The results are finally proven in Section \ref{proofs}.\\

\noindent
{\it Acknowledgement}. 
The authors would like to  thank Roger Tribe for many fruitful discussions.
\section{Notation}\label{notation}
We use the notation $\partial_i$ for the operation of taking the $i$th
partial derivative, $i=1,\dots,d$, of a function $f(x_1,\dots,x_d)$
and we write $\partial_i^m$ for iterating this operation $m$ times,
that is, $\partial_i^m=\partial_i\partial_i^{m-1}$ for $m\ge 1$ where
$\partial_i^0$ is defined to be the identity map. But if the function
$f$ only depends on a space variable $x\in\NR$ 
and a time variable $t\ge 0$ and there is no ambiguity about the
nature of the involved variables
then we will also write $\partial_x$ and $\partial_t$ for the corresponding
partial derivatives.

The heat kernel
$$g(y,s\,;x,t)\,\stackrel{\mbox{\tiny def}}{=}\,
\frac{1}{\sqrt{4\pi(t-s)}}\,
\exp\{\frac{-(x-y)^2}{4(t-s)}\}\,
{\bf 1}_{(s,\infty)}(t)$$
is considered a function
$$g:\left[\rule{0pt}{10pt}\right.
\NR\times\NR
\left.\rule{0pt}{10pt}\right]
\times
\left[\rule{0pt}{10pt}\right.
\NR\times\NR
\left.\rule{0pt}{10pt}\right]
\to\NR.$$
We write $g$ as an inhomogeneous transition kernel in order to emphasise that
the method works for more general PDE problems than the heat equation.
However for some explicit calculations we are going to apply the time\,-\,homogeneous structure of $g$
and then we also use the notation
$$g_y^{x_0}(t)\,\stackrel{\mbox{\tiny def}}{=}\,
\frac{1}{\sqrt{4\pi t}}\exp\{-\frac{(x_0-y)^2}{4t}\}\,\ind_{(0,\infty)}(t),
\quad t\in\NR,$$
for given $x_0,y\in\NR$.

We use $f_1\ast f_2$ to denote the convolution of functions $f_1,f_2:\NR\to\NR$
and write $\widehat{f_i}$ for their Laplace transform
$$\widehat{f_i}(\nu)\,\stackrel{\mbox{\tiny def}}{=}\,
\int_0^\infty f_i(t)e^{-\nu t}\,\dd t,\quad\nu>0,\;i=1,2.$$
Note that 
$$\widehat{(f_1\ind_{(0,\infty)})\ast(f_2\ind_{(0,\infty)})}\,=\,
\widehat{f_1}\,\widehat{f_2}.$$
Furthermore, if $l$ is a function on $(0,\infty)$ or $[0,\infty)$
then we denote by $l^a$ its antisymmetric extension
$$l^a:\NR\to\NR
\quad\mbox{such that $l^a(0)=0$ and}\quad
l^a(t)\,=\left\{\begin{array}{rcl}
l(t)&:&t\in(0,\infty);\\-l(-t)&:&t\in(-\infty,0).
\end{array}\right.$$
Note that 
$\|l^a\|_{L^2(\NR)}=\sqrt{2}\|l\|_{L^2([0,\infty))}$,
$\|l^a\|_{L^1(\NR)}={2}\|l\|_{L^1([0,\infty))}$
and that
\begin{eqnarray}\label{Young}
\|\frac{1}{\sqrt{|\cdot|}}\ast l^a\|_{L^p(\NR)}
&\le&
\|\frac{\ind_{[-1,1]}(\cdot)}{\sqrt{|\cdot|}}\ast l^a\|_{L^p(\NR)}
+
\|\frac{\ind_{\NR\setminus[-1,1]}(\cdot)}{\sqrt{|\cdot|}}\ast l^a\|_{L^p(\NR)}
\nonumber\\
&\le&
c_p(\|l^a\|_{L^2(\NR)}+\|l^a\|_{L^1(\NR)})
\end{eqnarray}
for each $p>2$ by Young's inequality.

The following domains
$${\cal Q}_+\,=\,\NR\times(0,\infty)
\quad\mbox{and}\quad
{\cal Q}_+^y\,=\,(y,\infty)\times(0,\infty),\;y\in\NR,$$
will be frequently used.

The symbol ${\mathscr D}$ is reserved for 
$C_c^\infty\hspace{-3pt}\left(\rule{0pt}{10pt}\right.
        \hspace{-3pt}(0,\infty)\hspace{-3pt}\left.\rule{0pt}{10pt}\right)$
the space of smooth functions on $(0,\infty)$ with compact support
and ${\mathscr D}^a\stackrel{\mbox{\tiny def}}{=}\{h^a:h\in{\mathscr D}\}$.

$\langle\cdot\,;\cdot\rangle$ denotes the scalar product in $L^2([0,\infty))$ and,
whenever the dual pairing between a topological vector space and its dual
is an extension of the scalar product in $L^2([0,\infty))$,
this dual pairing is also denoted by $\langle\cdot\,;\cdot\rangle$.
\section{Results}\label{results}
The stochastic partial differential equation of our interest
formally reads
\begin{equation}\label{our SPDE}
\partial_t u\,=\,\partial_x^2 u
\,+\,\partial_x\partial_t B\quad\mbox{in}\;{\cal Q}_+
\quad\mbox{subject to}\quad
\mbox{$\lim_{t\downarrow 0}$}\,u(\cdot,t)\,=\,0
\end{equation}
where $B=\{B_{xt};\,x\in\NR,t\ge 0\}$ is a Brownian sheet
given on a complete probability space $(\Omega,{\cal F},\pp)$.
Note that the transition kernel $g$ introduced in Section \ref{notation}
gives the Green's function associated with this Cauchy problem.
It is well-known that the random field
\begin{equation}\label{greenrep}
U(x,t)\,\stackrel{\mbox{\tiny def}}{=}\,
\dB\,g(y,s\,;x,t),
\quad(x,t)\in{\cal Q}_+,
\end{equation}
is the unique (weak) solution to (\ref{our SPDE})
where the integral against
$B(\dd y,\dd s)$ is understood as an It\^o-type integral 
against a process indexed by two parameters.
We always mean by $U$ the version which can be 
continuously extended to the closure of ${\cal Q}_+$ -- see \cite{W1986} for a good account
on the underlying theory of SPDEs.

Due to the parabolic nature of (\ref{our SPDE}), the process
$\{U(\cdot,t);\,t\ge 0\}$ taking values in the space of continuous functions
is a strong Markov process with zero initial condition in the usual sense.
But we are after the Markovian behaviour
of $U(x,\cdot)$ as a process indexed by $x\ge x_0$.
The method is to construct a (infinite dimensional) stochastic differential equation
wich is solved by the pair $(U(x,\cdot),\partial_1 U(x,\cdot))$ and then to prove that
the solution of this stochastic differential equation is Markovian 
in the ususal sense. 
\begin{rem}\rm
It turns out that this Markov process is homogeneous and stationary
in the case of (\ref{our SPDE}), that is, in the case of zero initial condition.
If the initial condition is not zero but a function $b_0$ then 
the underlying solution can be written as $U(x,t)+\int_\NR b_0(y)g(y,0\,;x,t)\,\dd y$
if $b_0$ has enough regularity. 
Adding this deterministic integral gives an inhomogeneous Markov process instead
without any extra proof.
\end{rem}

The initial idea is to rewrite (\ref{our SPDE}) as
\begin{equation}\label{system}
\left.
\begin{array}{rcl}
\displaystyle
\partial_x u&=&\displaystyle v\\
\displaystyle
\partial_x v&=&\displaystyle
\partial_t u-\partial_x\partial_t B
\end{array}
\right\}
\end{equation}
and to understand this system as an equivalent formulation of the Dirichlet problem
\begin{equation}\label{dirichlet}
(\partial_x^2-\partial_t)\,u
\,=\,-\partial_x\partial_t B
\quad\mbox{in}\quad
{\cal Q}_+^{x_0}\,=\,
(x_0,\infty)\times(0,\infty)
\end{equation}
subject to a given continuous function on the boundary $\partial{\cal Q}_+^{x_0}$.
So we are only interested in solutions $(u,v)$ of (\ref{system}) with respect
to the domain ${\cal Q}_+^{x_0}$ where $u(x,t)$ can be extended to a 
continuous function on $\overline{{\cal Q}_+^{x_0}}$ satisfying
\begin{equation}\label{boundary}
u(x_0,\cdot)\,=\,b^{x_0}
\quad\mbox{and}\quad
u(\cdot,0)\,=\,b_0
\end{equation}
for given (maybe random) continuous functions 
$$b^{x_0}:[0,\infty)\to\NR
\quad\mbox{and}\quad 
b_{0}:[x_0,\infty)\to\NR
 \quad\mbox{such that}\quad
b^{x_0}(0)=b_0(x_0).$$
\begin{rem}\rm\label{sol dirichlet}
\begin{itemize}\item[(i)]
If a continuous function $u$ on ${\cal Q}_+^{x_0}$ satisfies (\ref{dirichlet})
in the weak sense of
$$\int\hspace{-10pt}\int_{{\cal Q}_+^{x_0}}\hspace{-0pt}
u(x,t)\,(\partial_1^2+\partial_2)f(x,t)\,\dd x\dd t
\,\stackrel{\mbox{\tiny\rm a.s.}}{=}\,
-\int\hspace{-10pt}\int_{{\cal Q}_+}\hspace{-10pt}B(\dd x,\dd t)
\,f(x,t)$$
for all test functions $f\in C_c^\infty({\cal Q}_+^{x_0})$
then the pair $(u,v)$ where $v$ is the generalised function given by
$v(f)=-\int\hspace{-8pt}\int_{{\cal Q}_+^{x_0}}\hspace{-0pt}
u(x,t)\,\partial_1 f(x,t)\,\dd x\dd t,\,f\in C_c^\infty({\cal Q}_+^{x_0})$,
satisfies (\ref{system}) in the corresponding weak sense and vice versa.
\item[(ii)]
If $b^{x_0}$ is an arbitrary 
continuous function on $[0,\infty)$
and $b_0$ is a continuous function on $[x_0,\infty)$ of polynomial growth 
such that $b^{x_0}(0)=b_0(x_0)$ then
\begin{eqnarray*}
u(x,t)&=&
\int_0^\infty b^{x_0}(s)\,2\partial_1 g(x_0,s\,;x,t)\,\dd s\\
&+&\rule{0pt}{20pt}
\int_{x_0}^\infty b_0(y)\,[g(y,0\,;x,t)-g(2x_0-y,0\,;x,t)]\,\dd y\\
&+&\rule{0pt}{20pt}
\int\hspace{-10pt}\int_{{\cal Q}_+^{x_0}}\hspace{-0pt}B(\dd y,\dd s)\,
[g(y,s\,;x,t)-g(2x_0-y,s\,;x,t)]
\end{eqnarray*}
is the unique (weak) solution of (\ref{dirichlet}),(\ref{boundary}).
Using the Green's function associated with this Dirichlet problem,
the above existence/uniqueness result is standard -- see \cite{W1986} for example.
The polynomial growth condition on $b_0$ is not optimal but sufficient for our purpose.
Both, $b^{x_0}$ and $b_0$, can of course be taken to be ${\cal F}$\,-\,measurable.
\end{itemize}
\end{rem}

Let us return to the solution $U$ of (\ref{our SPDE}) given by (\ref{greenrep}).
Of course, $U$ is the unique weak solution of (\ref{dirichlet}),(\ref{boundary})
subject to $b^{x_0}=U(x_0,\cdot)$ and $b_0=0$ hence, by Remark \ref{sol dirichlet},
the pair $(U,\partial_1 U)$ solves (\ref{system}) at least in the corresponding
weak sense.

But this is not enough to justify why $(U(x,\cdot),\partial_1 U(x,\cdot))$ 
should be a Markov process indexed by $x\ge x_0$.
First one needs a meaning
of $(U(x,\cdot),\partial_1 U(x,\cdot))$ as a process indexed by $x\ge x_0$
which boils down to finding an appropriate state space, $E$, for 
the random variables $(U(x,\cdot),\partial_1 U(x,\cdot)),\,x\ge x_0$.
Second, a well-posed martingale problem needs to be associated 
with the system (\ref{system}).

To start with finding the right state space, observe that
$$\int_\NR\int_0^\infty\left(
\left[\rule{0pt}{10pt}\right.\hspace{-3pt}
\int_0^\infty\hspace{-10pt}\,g(y,s\,;x,t){h}(t)\,\dd t
\hspace{-3pt}\left.\rule{0pt}{10pt}\right]^2
\,+\,
\left[\rule{0pt}{10pt}\right.\hspace{-3pt}
\int_0^\infty\hspace{-10pt}\partial_3 g(y,s\,;x,t){h}(t)\,\dd t
\hspace{-3pt}\left.\rule{0pt}{10pt}\right]^2
\right)\dd s\dd y\,<\,\infty$$
hence the stochastic integrals
\begin{equation}\label{defi U}
U(x,h)\,\stackrel{\mbox{\tiny def}}{=}\,
\dB
\left[\rule{0pt}{10pt}\right.\hspace{-3pt}
\int_0^\infty\hspace{-10pt}\,g(y,s\,;x,t){h}(t)\,\dd t
\hspace{-3pt}\left.\rule{0pt}{10pt}\right]
\end{equation}
and
\begin{equation}\label{defi partial U}
\partial_1 U(x,h)\,\stackrel{\mbox{\tiny def}}{=}\,
\dB
\left[\rule{0pt}{10pt}\right.\hspace{-3pt}
\int_0^\infty\hspace{-10pt}\partial_3 g(y,s\,;x,t){h}(t)\,\dd t
\hspace{-3pt}\left.\rule{0pt}{10pt}\right]
\end{equation}
are well-defined
for all $x\ge x_0$ and 
$h\in{\mathscr D}=C_c^\infty\hspace{-3pt}\left(\rule{0pt}{10pt}\right.
        \hspace{-3pt}(0,\infty)\hspace{-3pt}\left.\rule{0pt}{10pt}\right)$.
Since the notation
$\partial_1^0 U(x,h)$ and $\partial_1^1 U(x,h)$
can be used for
$U(x,h)$ and $\partial_1 U(x,h)$,
respectively,
we have defined 
a centred Gaussian process $\partial_1^i U(x,h)$ 
indexed by $(i,x,h)\in\{0,1\}\times[x_0,\infty)\times{\mathscr D}$.
\begin{prop}\label{stateSpace}
\begin{itemize}\item[(i)]
The process $\{\partial_1^i U(x,h);\,(i,x,h)\in\{0,1\}\times[x_0,\infty)\times{\mathscr D}\}$
solves the system (\ref{system}) in the sense of
\begin{eqnarray*}
U(x,h)&\stackrel{\mbox{\tiny\rm a.s.}}{=}&
U(x_0,h)\,+
\int_{x_0}^x\partial_1 U(y,h)\,\dd y\\
\partial_1 U(x,h)&\stackrel{\mbox{\tiny\rm a.s.}}{=}&
\partial_1 U(x_0,h)\,-
\int_{x_0}^x U(y,h')\,\dd y
\,-\,
\dB\,(\ind_{(x_0,x]}\otimes h)(y,s)
\end{eqnarray*}
for all $(x,h)\in[x_0,\infty)\times{\mathscr D}$.
\item[(ii)]
For fixed $x\ge x_0$,
the processes $\{U(x,h);\,h\in{\mathscr D}\}$ and $\{\partial_1 U(x,h);\,h\in{\mathscr D}\}$
are independent centred Gaussian processes with covariances
$$\ee\,U(x,h_1)U(x,h_2)\,=\,
\langle\,h_1\,;\frac{-\sqrt{|\cdot|}}{\sqrt{4\pi}}\ast h_2^a\,\rangle$$
and
$$\ee\,\partial_1 U(x,h_1)\partial_1 U(x,h_2)\,=\,
\langle\,h_1\,;\frac{1}{2\sqrt{4\pi|\cdot|}}\ast h_2^a\,\rangle$$
respectively.
\item[(iii)]
For fixed $x\ge x_0$,
the process $\{U(x,h);\,h\in{\mathscr D}\}$
has a version taking values in
$$E_1\,\stackrel{\mbox{\tiny def}}{=}\,\{u\in(C_{0,\alpha})':
\mbox{$u\in C([0,\infty))$ such that $u(0)=0$}\}$$
for some $\alpha>3/2$ where 
$$C_{0,\alpha}\,\stackrel{\mbox{\tiny def}}{=}\,
\{h\in C([0,\infty)):h(t)(1+t)^{\alpha}\to 0,\,t\to\infty\},$$
is equipped with the norm 
$\|h\|_{0,\alpha}\stackrel{\mbox{\tiny def}}{=}\sup_{t\ge 0}|h(t)(1+t)^{\alpha}|$
and $(C_{0,\alpha})^\prime$ denotes the topological dual 
of the Banach space $(C_{0,\alpha}\,,\|\cdot\|_{0,\alpha})$.
Equip $E_1$ with the norm of $(C_{0,\alpha})'$.
\item[(iv)]
For fixed $x\ge x_0$, 
the process $\{\partial_1 U(x,h);\,h\in{\mathscr D}\}$ has a version taking values in
$$E_2\,\stackrel{\mbox{\tiny def}}{=}\,
(H_{\mathfrak{w},\beta}^a)'
\quad\mbox{for some $\beta>1/4$}$$
where
$$H_{\mathfrak{w},\beta}^a\,=\,\{l\in L^2([0,\infty)):\mathfrak{w}l^a\in H_\beta\}.$$
Here $H_\beta$ stands for the Sobolev space of functions $f\in L^2(\NR)$
whose Fourier transform $f^F$ satisfies
$\|(1+|\cdot|^2)^{\frac{\beta}{2}}f^F\|_{L^2(\NR)}\,<\,\infty$
and $\mathfrak{w}$ is a smooth weight function such that, for some $\varepsilon>0$,
$\mathfrak{w}\ge 1+|\cdot|^{\frac{1}{2}+\varepsilon}$ 
but
$\mathfrak{w}=1+|\cdot|^{\frac{1}{2}+\varepsilon}$ outside a neighbourhood of zero.
\item[(v)]
The family of random variables $\{(U(x,\cdot),\partial_1 U(x,\cdot));\,x\ge x_0\}$ 
is a stationary process taking values in $E=E_1\times E_2$
which has an ${\cal F}\otimes{\mathscr B}([x_0,\infty))$\,-\,measurable version
such that
\begin{equation}\label{bound on norms}
\ee\int_{x_0}^x\left(\,
\|U(y,\cdot)\|_{E_1}^2+\|\partial_1 U(y,\cdot)\|_{E_2}^2
\,\right)\,\dd y
\,<\,\infty
\end{equation}
for all $x>x_0$. 
Furthermore, the process $\{U(x,\cdot);\,x\ge x_0\}$ taking values in $E_1$
has a version such that $(x,t)\mapsto U(x,t)$ is continuous
on the closure of ${\cal Q}_+^{x_0}$.
\end{itemize}
\end {prop}
\begin{rem}\label{why H minus two}\rm
\begin{itemize}\item[(i)]
The bound $3/2$ for the parameter $\alpha$ defining the state space $E_1$
is sharp in the following sense: for $\alpha\le 3/2$ one cannot apply
Lemma \ref{fubini} below in the proof.
\item[(ii)]
Denote by $C_2$ the covariance operator
$C_2 h\,\stackrel{\mbox{\tiny def}}{=}\,
\frac{1}{2\sqrt{4\pi|\cdot|}}\ast h^a$ associated with $\partial_1 U(x,\cdot)$ by item (ii) above.
Of course, $C_2 h=const(-\partial_t^2)^{-\frac{1}{4}}\,h^a$ (see \cite{S1970} for example), 
and the parameter $\beta$ defining the space $E_2$ was chosen just big enough to ensure that 
$C_2^{1/2}:L^2([0,\infty))\to E_2$ is a Hilbert-Schmidt operator which,
by Sazonov's theorem, is needed for a meaningful state space of a Gaussian measure.
The choice of a weighted (Sobolev) space is due to the `infinite-volume' in $t$-direction. 
Finding the right space $E_2$ in the case of other SPDEs might be more complicated.
\end{itemize}
\end{rem}

The proof of the above proposition uses the following technical lemma.
Recall that $B=\{B_{ys};\,y\in\NR,s\ge 0\}$ is a Brownian sheet
on a given complete probability space $(\Omega,{\cal F},\pp)$. 
We assign to $B$ a family of $\sigma$-algebras
$$\label{brownian sigma field}
{\cal F}_{\!A}\,=\,\sigma(\{B_{ys}:y\in A,s>0\})
\vee{\cal N}_\pp,\quad A\subseteq\NR,$$
where ${\cal N}_\pp$ is the collection of all null sets in ${\cal F}$.
Note that this makes ${\cal F}_{(-\infty,x]},\,x\in\NR$, a right-continuous filtration.
\begin{lemma}\label{fubini}
{\rm(special case of \cite[Th.\ 2.6]{W1986})}.
Let $\phi\in L^1(I)$ 
where $I\subseteq\NR$ is a measurable index set
and let $f:\Omega\times{\cal Q}_+\times I\to\NR$
be an ${\cal F}\otimes{\mathscr B}({\cal Q}_+)\otimes{\mathscr B}(I)$\,-\,measurable
function such that, for each $(y,\zeta)\in\NR\times I$, 
the mapping $(\omega,s)\mapsto f(\omega,y,s,\zeta)$ is 
${\cal F}_{(-\infty,y]}\otimes{\mathscr B}((0,\infty))$\,-\,measurable.
\begin{itemize}
\item[(i)]
If
$\ee\int_\NR\int_0^\infty[f(y,s,\zeta)]^2\,\dd s\dd y<\infty$
for all $\zeta\in I$
then there is an ${\cal F}\otimes{\mathscr B}(I)$\,-\,measurable
version of the process 
$\{\iint_{{\cal Q}_+}B(\dd y,\dd s)\,f(y,s,\zeta);\,\zeta\in I\}$.
\item[(ii)]
If in addition
$$\ee\int_\NR\int_0^\infty\int_I
[f(y,s,\zeta)]^2\,|{\phi}(\zeta)|\,\dd\zeta\,\dd s\dd y\,<\,\infty$$
then the integrals below exist and satisfy
$$\int_I
\left[\rule{0pt}{10pt}\right.\hspace{-3pt}
\dB\,f(y,s,\zeta)
\hspace{-3pt}\left.\rule{0pt}{10pt}\right]
\phi(\zeta)\,\dd\zeta
\,\stackrel{\mbox{\tiny\rm a.s.}}{=}\,
\dB\left[\rule{0pt}{10pt}\right.\hspace{-3pt}
\int_I\hspace{-0pt}f(y,s,\zeta)\,\phi(\zeta)\,\dd\zeta
\hspace{-3pt}\left.\rule{0pt}{10pt}\right].$$
\end{itemize}
\end{lemma}

Now we introduce the random variables
$$\label{defi wipro}
W_z(l)\,\stackrel{\mbox{\tiny def}}{=}\,
\dB\,(\ind_{(x_0,x_0+z]}\otimes l)(y,s),
\quad(z,l)\in[0,\infty)\times L^2([0,\infty)),$$
such that the equations of Proposition \ref{stateSpace}(i) can be rewritten as
\begin{equation}\label{basic SDE}
\left. \begin{array}{rcl}
U(x,h)&\stackrel{\mbox{\tiny\rm a.s.}}{=}&
U(x_0,h)\,+
\int_{x_0}^x\partial_1 U(y,h)\,\dd y\\
\rule{0pt}{20pt}
\partial_1 U(x,h)&\stackrel{\mbox{\tiny\rm a.s.}}{=}&
\partial_1 U(x_0,h)\,-
\int_{x_0}^x U(y,h')\,\dd y
\,-\,
W_{x-x_0}(h)
\end{array}\right\}
\end{equation}
for all $(x,h)\in[x_0,\infty)\times{\mathscr D}$.
Note that, for fixed $l\in L^2([0,\infty))\setminus\{0\}$,
a version of the process $\{W_z(l)/\|l\|_{L^2([0,\infty))};\, z\ge 0\}$ is a standard  Wiener process
with respect to the filtration ${\cal F}_{(-\infty,x_0+z]},\,z\ge 0$,
and $\{W_z(0)=0;\,z\ge 0\}$ is a version for $l=0$.
These versions are used for all processes of type $\{aW_z(l);\,z\ge 0\}$
with fixed $(a,l)\in\NR\times L^2([0,\infty))$ in what follows.

So, if the process 
$\{(U(x,\cdot),\partial_1 U(x,\cdot));\,x\ge x_0\}$
were ${\cal F}_{(-\infty,x]}$\,-\,adapted
then one could try to establish the Markov property of this process via the  
martingale problem corresponding to the stochastic differential equation (\ref{basic SDE}).
But, from (\ref{defi U}) follows that
$$\ee[U(x,h)\,|\,{\cal F}_{(-\infty,x]}]
\,\stackrel{\mbox{\tiny\rm a.s.}}{=}\,
\int\hspace{-10pt}\int_{{\cal Q}_+\setminus{\cal Q}_+^{x}}\hspace{-0pt}B(\dd y,\dd s)
\left[\rule{0pt}{10pt}\right.\hspace{-3pt}
\int_0^\infty\hspace{-10pt}\,g(y,s\,;x,t){h}(t)\,\dd t
\hspace{-3pt}\left.\rule{0pt}{10pt}\right]
\,\not=\,U(x,h)$$
for any $(x,h)\in[x_0,\infty)\times{\mathscr D}$
thus $\{(U(x,\cdot),\partial_1 U(x,\cdot));\,x\ge x_0\}$ cannot be ${\cal F}_{(-\infty,x]}$\,-\,adapted.

The crucial observation is now that 
$\{(U(x,\cdot),\partial_1 U(x,\cdot));\,x\ge x_0\}$ is adapted with respect to the enlarged filtration
$$\tilde{\cal F}_x\,\stackrel{\mbox{\tiny def}}{=}\,
{\cal F}_{(-\infty,x]}\vee\sigma(U(x_0,\cdot)),\quad x\ge x_0.$$
The intuition behind this is of course that a unique solution to (\ref{basic SDE}) should be a functional
of the initial data $U(x_0,\cdot),\,\partial_1 U(x_0,\cdot)$
and the driving Wiener process. In our case this can easily be made precise
by approximating the derivative $h'$ in (\ref{basic SDE}) by a bounded operator
and showing that the $\tilde{\cal F}_x$\,-\,adapted solutions of the approximating systems converge to the
unique solution of (\ref{basic SDE}). To do so, we would use the connection
between (\ref{basic SDE}) and (\ref{dirichlet}) as explained in Remark \ref{sol dirichlet}.
The wanted convergence can then be verified in a straight forward way
using the Green's function given by Remark \ref{sol dirichlet}(ii).

As a consequence, $\{(U(x,\cdot),\partial_1 U(x,\cdot));\,x\ge x_0\}$ is at least adapted
with respect to the filtration
${\cal F}_{(-\infty,x]}\vee\sigma(U(x_0,\cdot),\partial_1 U(x_0,\cdot))$
but, by Remark \ref{partial_1 U adaptness} on page \pageref{partial_1 U adaptness},
we know that $\partial_1 U(x_0,\cdot)$ is
$\tilde{\cal F}_{x_0}$\,-\,measurable so that 
$\{(U(x,\cdot),\partial_1 U(x,\cdot));\,x\ge x_0\}$ is indeed $\tilde{\cal F}_x$\,-\,adapted.
Note that, in the case of other SPDEs, it can easily happen that one has to enlarge
${\cal F}_{(-\infty,x]}$ by initial conditions with respect to several
partial derivatives.

However, since $\{W_z(l);\,z\ge 0\}$ is not a martingale
with respect to the bigger filtration $\tilde{\cal F}_{x_0+z},\,z\ge 0$,
the equation (\ref{basic SDE}) cannot be associated with a martingale problem
in a straight forward way, yet.
One has to find a semimartingale decomposition 
of the process $\{W_z(l);\,z\ge 0\}$ with respect to
$\tilde{\cal F}_{x_0+z},\,z\ge 0$, 
and this problem is dealt with in the next proposition.

First we state the well-known martingale representation theorem
for  Brownian sheet: 
if ${\cal L}$ is an ${\cal F}_\NR$\,-\,measurable random variable in $L^2(\Omega)$
then there exists an ${\cal F}_{(-\infty,y]}$\,-\,adapted measurable process
$(\dot{\lambda}_{y\,\cdot})_{y\in\NR}$ in
$L^2(\Omega\times\NR\,;L^2([0,\infty))$ such that
$${\cal L}\,\stackrel{\mbox{\tiny\rm a.s.}}{=}\,\ee{\cal L}\,+\dB\,\dot{\lambda}_{ys}\,.$$
A good reference for this result in an even more general setting is \cite[Th.1.4]{N2006}.

As a consequence, any ${\cal F}_\NR$\,-\,measurable random variable $L$
taking values in a measurable space $E$
is associated with an additive  stochastic kernel $\dot{\lambda}_{ys}(F)$ 
indexed by bounded measurable functions 
$F:E\to\NR$ such that
$$F(L)\,\stackrel{\mbox{\tiny\rm a.s.}}{=}\,
\ee F(L)\,+\dB\,\dot{\lambda}_{ys}(F)\,.$$ 

In what follows let $E$ be a Souslin locally convex topological vector space 
and denote by $E'$ its topological dual. Introduce
\begin{eqnarray*} 
{\mathfrak F}C_b^\infty(D)
&\stackrel{\mbox{\tiny def}}{=}&
\left\{\rule{0pt}{13pt}\right.
\mbox{$F:E\to\NR$ such that $F(\phi)=f(h_1(\phi),\dots,h_n(\phi))$ for}\\
&&\rule{0pt}{20pt}\hspace{.4cm}
f\in C_b^\infty(\NR^n),\,h_i\in D,\,i=1,\dots,n,\,n\in\NN
\left.\rule{0pt}{13pt}\right\}
\end{eqnarray*} 
where $D\subseteq E'$ is supposed to separate the points of $E$.
Then we have both
$\sigma(L)\,=\,\sigma(\{F(L):F\in{\mathfrak F}C_b^\infty(D)\})$
and $\{F(L):F\in{\mathfrak F}C_b^\infty(D)\}$ is dense in 
$L^2(\Omega,\sigma(L),\pp)$ so that the kernel $\dot{\lambda}_{ys}(F)$
is fully described by $F\in{\mathfrak F}C_b^\infty(D)$.
\begin{prop}\label{new wipro}
Fix an ${\cal F}_\NR$\,-\,measurable random variable $L:\Omega\to E$
and $l\in L^2([0,\infty))$.
Assume that there exists a measurable function
$$\varrho_l:\Omega\times E\times[x_0,\infty)\to\NR$$
such that
\begin{itemize}
\item
$\varrho_l(\phi,y)$ is ${\cal F}_{(-\infty,y]}$\,-\,measurable
for each $\phi\in E$ and $y\ge x_0$\,;
\item
$\varrho_l(L,y)\in L^1(\Omega)$ for almost every $y\ge x_0$\,;
\item
the mapping $y\mapsto\varrho_l(L,y)$ is in $L^1([x_0,x])$ almost surely
for each $x>x_0$\,;
\item
for each $F\in{\mathfrak F}C_b^\infty(D)$
and almost every $y\ge x_0$
it holds that
\begin{equation}\label{how rho}
\int_0^\infty\dot{\lambda}_{ys}(F)\,l(s)\,\dd s\,\stackrel{\mbox{\tiny\rm a.s.}}{=}\,
\ee\left[\left.\rule{0pt}{11pt}
F(L)\varrho_l(L,y)\right|{\cal F}_{(-\infty,y]}
\right].
\end{equation}
\end{itemize}
If
$$\tilde{W}_z(l)\,\stackrel{\mbox{\tiny def}}{=}\,
W_z(l)\,-\int_{x_0}^{x_0+z}\varrho_l(L,y)\,\dd y,\quad z\ge 0,$$
then, for $l\not=0$, the process
$\{\tilde{W}_z(l)/\|l\|_{L^2([0,\infty))};\,z\ge 0\}$ is a standard Wiener process
with respect to the filtration ${\cal F}_{(-\infty,x_0+z]}\vee\sigma(L),\,z\ge 0$.
Moreover, if $\varrho_l$ with the above properties exists for
$l=l_1,l_2$ then $\varrho_{a_1 l_1+a_2 l_2}$ exists for each $a_1,a_2\in\NR$ and
\begin{equation}\label{lin}
\tilde{W}_z(a_1 l_1+a_2 l_2)\,\stackrel{\mbox{\tiny\rm a.s.}}{=}\,
a_1\tilde{W}_z(l_1)+a_2\tilde{W}_z(l_2)
\end{equation}
for each $z\ge 0$.
\end{prop}
\begin{rem}\rm\label{idea new wipro}
\begin{itemize}\item[(i)]
This proposition is a generalisation of Theorem 12.1 in \cite{Y1997}
which deals with the semimartingale decomposition of a Wiener process $\{W_t;\,t\ge 0\}$
if its natural filtration ${\cal F}^W_t,\,t\ge 0$, is enlarged 
by the information given by an ${\cal F}^W_\infty$\,-\,measurable random variable.
In our case, for fixed $l\in L^2([0,\infty))\setminus\{0\}$,
the Wiener process $\{W_z(l)/\|l\|_{L^2([0,\infty))};\,z\ge 0\}$ 
is already a Wiener process with respect to a filtration larger than its natural filtration,
that is ${\cal F}_{(-\infty,x_0+z]},\,z\ge 0$,
and this larger filtration is enlarged further.
But we have both there is a martingale representation theorem 
with respect to ${\cal F}_{(-\infty,x_0+z]},\,z\ge 0$,
\underline{and}
$\{W_z(l);\,z\ge 0\}$ can be represented as a stochastic integral against the
${\cal F}_{(-\infty,x_0+z]}$\,-\,integrator which is the Brownian sheet.
So the idea of proof is the same as in the proof of \cite[Th.12.1]{Y1997}
so that, in the Proof-Section,
we will only deal with the following two elements of the proof:
the part where the different type of martingale representation is used
and the linearity (\ref{lin}).
\item[(ii)]
The proposition immediately implies that
if $\varrho^{}_l$ and $\varrho'_l$ are two functions satisfying all properties 
stated in the above proposition then 
$$\pp\left(\rule{0pt}{11pt}\right.
\varrho^{}_l(L,y)=\varrho'_l(L,y)\;for\;almost\;every\;y\ge x_0
\left.\rule{0pt}{11pt}\right)=1$$
because the process 
$\int_{x_0}^{x_0+z}(\varrho^{}_l(L,y)-\varrho'_l(L,y))\,\dd y,\,z\ge 0$,
is a continuous martingale and must vanish therefore.
\end{itemize}
\end{rem}

In our case, the role of $L$ in the above proposition is played by $U(x_0,\cdot)$
hence, by Proposition \ref{stateSpace}(iii), 
the corresponding Souslin locally convex space is $E_1$.
We choose ${\mathscr D}$ to be the subset of $E'_1$ separating the points of $E_1$.
The next lemma identifies a class of $l\in L^2([0,\infty))$ such that
$\varrho_l$ with the properties stated in Proposition \ref{new wipro}
exists for $L=U(x_0,\cdot)$.
\begin{lemma}\label{identify rho}
For an arbitrary but fixed $\nu>0$ set
$$l_\nu\,\stackrel{\mbox{\tiny def}}{=}\,
\frac{1}{\sqrt{4\pi|\cdot|}}\ast(e^{-\nu\,\fatdot}\,)^a.$$
Then $l_\nu$ is a bounded continuous function in $L^2([0,\infty))$ 
satisfying $l_\nu(0)=0$ and
$$\varrho_{l_\nu}(\phi,y)\,=\,
\langle\phi-U(x_0,\cdot)_y\,;
\frac{2\sqrt{\nu}e^{-\nu\,\fatdot}}{e^{-\sqrt{\nu}(y-x_0)}}\rangle,
\quad\phi\in E_1,\,y>x_0,$$
where
$$U(x_0,h)_y\,\stackrel{\mbox{\tiny def}}{=}\,
\int\hspace{-10pt}\int_{{\cal Q}_+\setminus{\cal Q}_+^{y}}\hspace{-0pt}B(\dd y',\dd s')
\left[\rule{0pt}{10pt}\right.\hspace{-3pt}
\int_0^\infty\hspace{-10pt}\,g(y',s'\,;x_0,t){h}(t)\,\dd t
\hspace{-3pt}\left.\rule{0pt}{10pt}\right].$$
Furthermore 
\begin{equation}\label{as drift}
\varrho_{l_\nu}(U(x_0,\cdot),y)\,\stackrel{\mbox{\tiny\rm a.s.}}{=}\,
U(y,\sqrt{\nu}e^{-\nu\,\fatdot})\,+\,\partial_1 U(y,e^{-\nu\,\fatdot}),\quad y>x_0,
\end{equation}
which does not depend on $x_0$ anymore.
\end{lemma}

The above lemma suggests that $\varrho_{l_\nu}(U(x_0,\cdot),y)$ can be written as
a sum of operators acting on $U(y,\cdot)$ and $\partial_1 U(y,\cdot)$ respectively.
In what follows, we will reveal the explicit nature of such operators.

For an absolutely continuous function $h:[0,\infty)\to\NR$ satisfying 
$h'\in L^1([0,\infty))\cap L^2([0,\infty))$ define the functions
${\mathfrak A}_1 h$ and ${\mathfrak A}_2 h$ on $[0,\infty)$ by
$${\mathfrak A}_1 h(t)\,=\int_t^\infty\frac{-h'(t')\,\dd t'}{\sqrt{\pi(t'-t)}}
\quad\mbox{and}\quad
{\mathfrak A}_2 h(t)\,=\,
\left[\rule{0pt}{10pt}\right.
\frac{sgn(\cdot)}{\sqrt{\pi|\cdot|}}\ast(h^a)'
\left.\rule{0pt}{10pt}\right](t)$$
respectively.
Since 
$|{\mathfrak A}_1 h|\le
\left[\rule{0pt}{10pt}\right.
|\cdot|^{-\frac{1}{2}}\ast|(h^a)'|
\left.\rule{0pt}{10pt}\right]$,
${\mathfrak A}_1 h$ and ${\mathfrak A}_2 h$ are well-defined by (\ref{Young}).
\begin{prop}\label{new drift}
\begin{itemize}\item[(i)]
If $h\in{\mathscr D}$ then the function $({\mathfrak A}_2 h)^a$ is a $C^\infty$-\,function 
satisfying 
$\partial_t^k{\mathfrak A}_2 h\in C_{0,\alpha}$ for $0\le\alpha<k+3/2,\,k=0,1,2,\dots$, 
and ${\mathfrak A}_2 h\in H_{\mathfrak{w},\beta}^a$ for all $\beta\ge 0$
if the parameter $\varepsilon>0$ used to determine the weight function $\mathfrak{w}$
in Proposition \ref{stateSpace}(iv) is less than $1/2$.
Also, the function ${\mathfrak A}_1{\mathfrak A}_2 h$ is in $C_{0,\alpha}$ for $0\le\alpha<2$.
\item[(ii)]
For $h\in{\mathscr D}$, the process
$$\tilde{W}_z(h)\,\stackrel{\mbox{\tiny def}}{=}\,
W_z(h)\,-\int_{x_0}^{x_0+z}[
U(y,{\mathfrak A}_1{\mathfrak A}_2 h)\,+\,
\partial_1 U(y,{\mathfrak A}_2 h)]\,\dd y,\quad z\ge 0,$$
is well-defined and if $h\not=0$ then
$\{\tilde{W}_z(h)/\|h\|_{L^2([0,\infty))};\,z\ge 0\}$
is a standard Wiener process with respect to the filtration $\tilde{\cal F}_{x_0+z},\,z\ge 0$.
Moreover, it holds that
\begin{equation}\label{lin new}
\tilde{W}_z(a_1 h_1+a_2 h_2)\,\stackrel{\mbox{\tiny\rm a.s.}}{=}\,
a_1\tilde{W}_z(h_1)+a_2\tilde{W}_z(h_2)
\end{equation}
for each $z\ge 0,\,a_1,a_2\in\NR,\,h_1,h_2\in{\mathscr D}$.
\item[(iii)]
If $h\in{\mathscr D}$ then
$${\mathfrak A}_1{\mathfrak A}_2 h\,=\,
-h'+(-\partial_t^2)^{\frac{1}{2}}\,h^a
\quad\mbox{and}\quad
{\mathfrak A}_2 h\,=\,\sqrt{2}\,(-\partial_t^2)^{\frac{1}{4}}\,h^a$$
where, for $\beta\in\NR$, the fractional Laplacian $(-\partial_t^2)^{\frac{\beta}{2}}f$
of $f\in C_c^\infty(\NR)$ is defined by its Fourier transform
$((-\partial_t^2)^{\frac{\beta}{2}}f)^F=|\cdot|^\beta f^F$.
\end{itemize}
\end{prop}
\begin{rem}\rm\label{why A1 and A2}
The operators ${\mathfrak A}_1$ and ${\mathfrak A}_2$ were introduced
to simplify the proof of item (ii) of the above proposition.
Furthermore, if $h\in{\mathscr D}$ then
$(-\partial_t^2)^{\frac{1}{2}}\,h^a\in C_{0,\alpha}$ for $0\le\alpha<2$ 
by Proposition \ref{new drift}(i) because
$(-\partial_t^2)^{\frac{1}{2}}\,h^a=
{\mathfrak A}_1{\mathfrak A}_2 h+h'$
and $h'$ has compact support.
\end{rem}

So, in what follows, we will always assume that the parameter $\varepsilon>0$ 
used to determine the weight function $\mathfrak{w}$
in Proposition \ref{stateSpace}(iv) is less than $1/2$.
Then, recalling (\ref{basic SDE}), Proposition \ref{new drift} implies that 
$\{(U(x,\cdot),\partial_1 U(x,\cdot));\,x\ge x_0\}$  
satisfies the equation
{
\begin{equation}\label{new SDE}
\left. \begin{array}{rcl}
U(x,h)&\stackrel{\mbox{\tiny\rm a.s.}}{=}&\displaystyle
U(x_0,h)\,+
\int_{x_0}^x\partial_1 U(y,h)\,\dd y\\
\rule{0pt}{20pt}
\partial_1 U(x,h)&\stackrel{\mbox{\tiny\rm a.s.}}{=}&\displaystyle
\partial_1 U(x_0,h)\,-
\int_{x_0}^x 
[\,U(y,(-\partial_t^2)^{\frac{1}{2}}\,h^a)+\partial_1 U(y,\sqrt{2}(-\partial_t^2)^{\frac{1}{4}}\,h^a)\,]
\,\dd y\\
&&\rule{0pt}{20pt}\displaystyle\hspace{1.8cm}
\,-\,
\tilde{W}_{x-x_0}(h)
\end{array}\right\}
\end{equation}
}\noindent
for all $(x,h)\in[x_0,\infty)\times{\mathscr D}$. Because both,
$\{(U(x,\cdot),\partial_1 U(x,\cdot));\,x\ge x_0\}$ is $\tilde{\cal F}_x$\,-\,adapted
and
$\{\tilde{W}_{z}(h);\,z\ge 0\}$ is a martingale with respect to $\tilde{\cal F}_{x_0+z}$,
this stochastic differential equation can eventually be associated with a martingale problem.
But before we do so let us point out that (\ref{new SDE}), 
when seen as a family of stochastic differential equations indexed by $x_0$,
gives raise to a new SPDE in ${\cal Q}_+$.
\begin{theo}\label{new SPDE}
The unique weak solution $U$ to (\ref{our SPDE}) given by the continuous version
of the right-hand side of (\ref{greenrep}) on page \pageref{greenrep} satisfies
$$\partial_x^2\,U +[\,(-\partial_t^2)^{1/2}+\sqrt{2}\,\partial_x(-\partial_t^2)^{1/4}\,]\,U^a=
\partial_x\partial_t{\tilde B}$$
in the sense of
$$U\left(\rule{0pt}{11pt}\right.
\partial_x^2 f+[\,(-\partial_t^2)^{1/2}-\sqrt{2}\,\partial_x(-\partial_t^2)^{1/4}\,]\,f^a
\left.\rule{0pt}{11pt}\right)
\,=\,
\partial_x\partial_t\tilde{B}(f)
\quad\mbox{for all $f\in C_c^\infty({\cal Q}_+)$}\quad a.s.$$
where $U^a,f^a$ stand for the extensions of $U(x,t),f(x,t)$ to $(x,t)\in\NR^2$ 
which are antisymmetric in $t$
and $\tilde{B}=\{\tilde{B}_{xt};\,x\in\NR,t\ge 0\}$ is a Brownian sheet on $(\Omega,{\cal F},\pp)$.
\end{theo}
\begin{rem}\rm\label{new meaning U}
In this theorem, $U$ is considered a regular generalised function on $C_c^\infty({\cal Q}_+)$, that is,
$$U(f)\,=\,\int_\NR\int_0^\infty U(x,t)\,f(x,t)\,\dd t\dd x
\quad\mbox{for}\quad f\in C_c^\infty({\cal Q}_+).$$
However, the proof of Proposition \ref{new drift}(i) makes clear that,
for fixed $x\in\NR$ and $f\in C_c^\infty({\cal Q}_+)$,
the long-time behaviour of $(-\partial_t^2)^{1/4}f^a(x,\cdot)$
is in general not better than ${\cal O}(t^{-3/2}),\,t\to\infty$, while, by Proposition \ref{stateSpace}(iii),
$U(x,\cdot)$ is only in $(C_{0,\alpha})'$ for $\alpha>3/2$.
So the meaning of $U(\,\partial_x(-\partial_t^2)^{1/4}f^a)$ is based on an extension
of the regular generalised function $U$ which will be explained in the proof of the theorem.
\end{rem}

Coming back to the martingale problem associated with (\ref{new SDE}),
choose $D={\mathscr D}\times{\mathscr D}$, which is a subset 
of the topological dual of $E=E_1\times E_2$, and denote by {\bf A} the subset of
$C_b(E)\times C(E)$ whose elements $(F,G)$ are given by
$$\label{set A}
F\in{\mathfrak F}C_b^\infty(D)
\quad\mbox{such that}\quad
F(\phi_1,\phi_2)\,=\,
f(\langle\phi_1\,;h_1\rangle,\langle\phi_2\,;h_2\rangle,\dots,
\langle\phi_1\,;h_{2n-1}\rangle,\langle\phi_2\,;h_{2n}\rangle)$$
for some $f\in C_b^\infty(\NR^{2n}),\,h_i\in{\mathscr D},\,i=1,2,\dots,2n$, and
{\small
\begin{eqnarray*}
G(\phi_1,\phi_2)&=&
\sum_{i\;odd}
\partial_i f(\dots,\langle\phi_1\,;h_i\rangle,\dots)\,
\langle\phi_2\,;h_i\rangle\\
&-&
\sum_{i\;even}
\partial_i f(\dots,\langle\phi_2\,;h_i\rangle,\dots)\,
[\,\langle\phi_1\,;(-\partial_t^2)^{\frac{1}{2}}\,h_i^a\rangle
+\langle\phi_2\,;\sqrt{2}(-\partial_t^2)^{\frac{1}{4}}\,h_i^a\rangle\,]\\
&+&\renewcommand{\arraystretch}{0.7}
\frac{1}{2}
\sum_{i,j\;even}
\partial_i\partial_j f(\dots,\langle\phi_2\,;
\begin{array}{c}\scriptstyle h_i\\
\scriptstyle or\\
\scriptstyle h_j\end{array}
\rangle,\dots,
\langle\phi_2\,;
\begin{array}{c}\scriptstyle h_j\\
\scriptstyle or\\
\scriptstyle h_i\end{array}\rangle,\dots)\,
\langle h_i\,;h_j\rangle.
\renewcommand{\arraystretch}{1.0}
\end{eqnarray*}
}\noindent
This definition of the subset {\bf A} of course requires 
$(-\partial_t^2)^{\frac{1}{2}}\,h^a\in E'_1$ and 
$(-\partial_t^2)^{\frac{1}{4}}\,h^a\in E'_2$
for $h\in{\mathscr D}$
which follows from Proposition \ref{new drift}(i) and Remark \ref{why A1 and A2}.

Then, according to \cite[Chapter 3]{EK1986}, by a
{\it solution of the martingale problem for {\bf A} with respect to ${\mathscr F}_z$}
one would 
mean an ${\mathscr F}_z$\,-\,progressively measurable process $R=(R_z)_{z\ge 0}$
on a filtered probability space $(\Omega,{\mathscr F},\pp)$ 
taking values in $E$ such that
$$F(R_{z})-F(R_{0})\,-\int_{0}^{z}G(R_y)\,\dd y,\quad z\ge 0,$$
is a martingale with respect to the filtration ${\mathscr F}_z,\,z\ge 0$,
for all $(F,G)\in{\bf A}$. When an initial condition $\mu$ is specified, 
it is also said that the process $R$ is a 
{\it solution of the martingale problem for $({\bf A},\mu)$}.

Next we check whether 
$\{(U(x_0+z,\cdot),\partial_1 U(x_0+z,\cdot));\,z\ge 0\}$
is a solution of the martingale problem for our set {\bf A} with respect to
$\tilde{\cal F}_{x_0+z}$.

First, the process $\{(U(x_0+z,\cdot),\partial_1 U(x_0+z,\cdot));\,z\ge 0\}$
has an $\tilde{\cal F}_{x_0+z}$\,-\,progressively measurable version
because it is $\tilde{\cal F}_{x_0+z}$\,-\,adapted and,
by Proposition \ref{stateSpace}(v),
has an ${\cal F}\otimes{\mathscr B}([x_0,\infty))$\,-\,measurable version
taking values in the space $E=E_1\times E_2$.
This can be verified the same way the analogous statement for real-valued 
adapted measurable processes was verified in \cite{CD1965}.
Notice that, by construction, the filtration $\tilde{\cal F}_x$
inherits right-continuity from the filtration ${\cal F}_{(-\infty,x]}$
defined on page \pageref{brownian sigma field}.

Second, knowing that 
$\{(U(x,\cdot),\partial_1 U(x,\cdot));\,x\ge x_0\}$
satisfies (\ref{bound on norms}) and (\ref{new SDE}), 
an easy application of It\^o's formula to 
$R_z=(U(x_0+z,\cdot),\partial_1 U(x_0+z,\cdot))$ yields that
$$F(R_{z})-F(R_{0})\,-\int_{0}^{z}G(R_y)\,\dd y,\quad z\ge 0,$$ 
is indeed a martingale with respect to $\tilde{\cal F}_{x_0+z}$
for all $(F,G)\in {\bf A}$.

Now we hope that the well-posedness-condition 
$$\left.\mbox{
\parbox{13cm}{
for each probability measure $\mu$ on $(E,{\mathscr B}(E))$,
any two solutions $R,R'$ of the martingale problem for $({\bf A},\mu)$
with respect to ${\mathscr F}_z,\,{\mathscr F}'_z$,
have the same one-dimensional distributions, that is, for each $z>0$,
$$\pp(\{R_z\in\Gamma\})\,=\,\pp'(\{R'_z\in\Gamma\}),\quad\Gamma\in{\mathscr B}(E),$$
}}\right\}\eqno({\rm wp})\label{wp}$$
as stated in \cite[Thm.4.2]{EK1986} is enough to ensure that a solution
$(R_z)_{z\ge 0}$ of the martingale problem is strong Markov in the sense of:
\begin{defi}\rm\label{defi strong markov}
Let $E$ be a separable metric space and let $\mu$ be a probability measure
on $(E,{\mathscr B}(E))$. An ${\mathscr F}_z$\,-\,progressively measurable
process $(R_z)_{z\ge 0}$ on a filtered probability space $(\Omega,{\mathscr F},\pp)$
taking values in $E$ is said to be a
{\it strong Markov process with initial condition $\mu$} if
\begin{itemize}
\item[(i)]
$\pp(\{R_{0}\in\Gamma\})\,=\,\mu(\Gamma)$ for all $\Gamma\in{\mathscr B}(E)$;
\item[(ii)]
for any ${\mathscr F}_z$\,-\,stopping time $\xi\ge 0,\,y\ge 0$ and $\Gamma\in{\mathscr B}(E)$,
$$\pp[\{R_{\xi+y}\in\Gamma\}\,|\,{\mathscr F}_\xi\,]
\,\stackrel{\mbox{\tiny\rm a.s.}}{=}\,
\pp[\{R_{\xi+y}\in\Gamma\}\,|\,\sigma(R_\xi)]
\quad\mbox{on $\{\xi<\infty\}.$}$$
\end{itemize}
\end{defi}

To show the strong Markov property of 
$R_z=(U(x_0+z,\cdot),\partial_1 U(x_0+z,\cdot))$
with initial condition
$\pp\circ(U(x_0,\cdot),\partial_1 U(x_0,\cdot))^{-1}$
we want to apply Thm.4.2(b) in \cite{EK1986}.
But the conclusion of this theorem is stated under the extra conditions that
${\bf A}\subseteq C_b(E)\times C_b(E)$ and that
$(R_z)_{z\ge 0}$ has a right-continuous version taking values in $E$.
\begin{rem}\rm\label{improve EK theorem}
\begin{itemize}\item[(i)]
Our set ${\bf A}$ defining the martingale problem is not a subset of $C_b(E)\times C_b(E)$
but of $C_b(E)\times C(E)$ only. However, in the general situation of \cite[Thm.4.2(b)]{EK1986},
the boundedness of $F$ and $G$ is the natural condition to ensure that
$|F(R_{z})-F(R_{0})\,-\int_{0}^{z}G(R_y)\,\dd y|$ has finite expectation for each $z\ge 0$.
In our specific situation, if $R_z=(U(x_0+z,\cdot),\partial_1 U(x_0+z,\cdot))$ then
$$\ee\left|\rule{0pt}{10pt}\right.\!
F(R_{z})-F(R_{0})\,-\int_{0}^{z}G(R_y)\,\dd y
\!\left.\rule{0pt}{10pt}\right|\,<\,\infty$$
for given $(F,G)\in{\bf A}$ because of (\ref{bound on norms}) and $F\in{\mathfrak F}C_b^\infty(D)$.
It turns out that Thm.4.2(b) in \cite{EK1986} remains valid
when adding a condition of type (\ref{bound on norms}) to the definition of 
the martingale problem--see Definition \ref{our mart problem}(iii) and Remark \ref{extra}(ii) below.
\item[(ii)]
Taking another look at the proof of Thm.4.2(b) in \cite{EK1986}
reveals that the right-continuous version of the solution
is only needed for
$$F(R_{z})-F(R_{0})\,-\int_{0}^{z}G(R_y)\,\dd y,\quad z\ge 0,$$
to be a right-continuous martingale
when $(F,G)\in{\bf A}$ in order to be able to apply Doob's optional sampling theorem.
So, it is already enough to require right-continuity of
$z\mapsto F(U(x_0+z,\cdot),\partial_1 U(x_0+z,\cdot))$ 
for all $F\in{\mathfrak F}C_b^\infty(D)$ to make the theorem
work in our case.
\end{itemize}
\end{rem}

It is easy to realise that there is a version of the process
$\{(U(x,\cdot),\partial_1 U(x,\cdot));\,x\ge x_0\}$
such that
$z\mapsto F(U(x_0+z,\cdot),\partial_1 U(x_0+z,\cdot))$
is continuous for all $F\in{\mathfrak F}C_b^\infty(D)$.
First, it is well-known (see \cite{W1986} for example)
that the process 
$\{W_z(l);\,(z,l)\in[0,\infty)\times L^2([0,\infty))\}$
defined on page \pageref{defi wipro} has a version
such that $\{W_z(\cdot);\,z\ge 0\}$ is a continuous ${\mathscr D}'$\,-\,valued process.
Second, using the above ${\mathscr D}'$\,-\,valued version 
of $\{W_z(\cdot);\,z\ge 0\}$ and the
${\cal F}\otimes{\mathscr B}([x_0,\infty))$\,-\,measurable version
of $\{(U(x,\cdot),\partial_1 U(x,\cdot));\,x\ge x_0\}$
as stated in Proposition \ref{stateSpace}(v),
one can construct from (\ref{basic SDE}) a version of
$\{(U(x,\cdot),\partial_1 U(x,\cdot));\,x\ge x_0\}$ such that
$$\pp\left(\rule{0pt}{11pt}\right.
x\mapsto U(x,h)\;\&\;x\mapsto\partial_1 U(x,h)\;
\mbox{are continuous for all $h\in{\mathscr D}$}
\left.\rule{0pt}{11pt}\right)
\,=\,1.
$$

So, by Thm.4.2(b) in \cite{EK1986} and Remark \ref{improve EK theorem},
the strong Markov property of our process
$\{(U(x_0+z,\cdot),\partial_1 U(x_0+z,\cdot));\,z\ge 0\}$
becomes a direct implication of the well-posedness-condition (wp).
But, for showing the uniqueness wanted in (wp),
we need to work with a more restrictive  martingale problem than
Ethier/Kurtz in \cite{EK1986}.
Recall the set ${\bf A}\subseteq C_b(E)\times C(E)$ introduced on page \pageref{set A}.
\begin{defi}\rm\label{our mart problem}
An ${\mathscr F}_z$\,-\,progressively measurable process $\{(u_z,v_z);\,z\ge 0\}$
on a filtered probability space $(\Omega,{\mathscr F},\pp)$ 
taking values in $E=E_1\times E_2$ is called a
{\it solution of the martingale problem for ${\bf A}$ with respect to ${\mathscr F}_z$}
iff
\begin{itemize}
\item[(i)]
the mappings $z\mapsto u_z(h)$ and $z\mapsto v_z(h)$ are continuous for all $h\in{\mathscr D}$;
\item[(ii)]
the map $(z,t)\mapsto u_z(t)$ is continuous on the closure of ${\cal Q}_+^0$\,;
\item[(iii)]
$\ee\int_0^z(\|u_y\|_{E_1}+\|v_y\|_{E_2})\,\dd y\,<\,\infty$ for all $z>0$;
\item[(iv)]
$\{F(u_z,v_z)-F(u_0,v_0)-\int_0^z G(u_y,v_y)\,\dd y;\,z\ge 0\}$
is a martingale with respect to the filtration ${\mathscr F}_z,\,z\ge 0$,
for all $(F,G)\in{\bf A}$.
\end{itemize}
\end{defi}
\begin{rem}\rm\label{extra}
\begin{itemize}\item[(i)]
We also use the phrase `solution of the martingale problem for $({\bf A},\mu)$'
when a specific initial condition $\mu$ is emphasised as in (wp).
\item[(ii)]
We claim that Thm.4.2(b) in \cite{EK1986} remains valid
with respect to our more restrictive definition of the martingale problem
when being applied to show the strong Markov property of 
$\{(U(x_0+z,\cdot)$, $\partial_1 U(x_0+z,\cdot));\,z\ge 0\}$.
A quick glance at the proof of this theorem shows that
one only has to pay attention to
the property (iii) of our Definition \ref{our mart problem}
and this will be done in the next item of this remark.
\item[(iii)]
Adapting the proof of Thm.4.2(b) in \cite{EK1986} to our setup,
fix a finite $\tilde{\cal F}_{x_0+z}$\,-\,stopping time $\xi$,
choose $\Xi\in\tilde{\cal F}_{x_0+\xi}$ such that $\pp(\Xi)>0$
and introduce
$$u_y\,\stackrel{\mbox{\tiny def}}{=}\,U(x_0+\xi+y,\cdot),\quad
v_y\,\stackrel{\mbox{\tiny def}}{=}\,\partial_1 U(x_0+\xi+y,\cdot)\quad
\mbox{for all $y\ge 0$}$$
and
$$\pp_1(\Gamma)\,\stackrel{\mbox{\tiny def}}{=}\,
\frac{\ee\ind_\Xi\pp[\Gamma|\tilde{\cal F}_{x_0+\xi}]}{\pp(\Xi)}\,,
\quad
\pp_2(\Gamma)\,\stackrel{\mbox{\tiny def}}{=}\,
\frac{\ee\ind_\Xi\pp[\Gamma|\sigma(U(x_0+\xi,\cdot),\partial_1 U(x_0+\xi,\cdot))]}{\pp(\Xi)}$$
for all $\Gamma\in{\cal F}$.
The task is to show the property in Definition \ref{our mart problem}(iii)
if $\ee$ is replaced by the expectation operators $\ee_1$ and $\ee_2$
given by the measures $\pp_1$ and $\pp_2$, respectively.
But, for fixed $z>0$, we obtain that
$$\ee_i\int_0^z(\|u_y\|_{E_1}+\|v_y\|_{E_2})\,\dd y\,\le\,
\frac{1}{\pp(\Xi)}\;
\ee\int_{x_0}^{x_0+\xi+z}(\|U(y,\cdot)\|_{E_1}+\|\partial_1 U(y,\cdot)\|_{E_2})\,\dd y$$
for $i=1,2$ where, by Proposition \ref{stateSpace}(v), the last term is finite
if the stopping time $\xi$ is bounded. 
And it is sufficient to check the strong Markov property
for bounded stopping times only--we refer to Problem 2.6.9 in \cite{KS1991} for example.
\end{itemize}
\end{rem}

After this preparation, the key part
of the proof of the below theorem consists in verifying the
well-posedness-condition (wp) on page \pageref{wp} for our martingale problem.
Recall the spaces $E_1,E_2$ defined in Proposition \ref{stateSpace}
and assume that the parameter $\varepsilon>0$ 
used to define $E_2$ is less than $1/2$.

\begin{theo}\label{strMtheo}
The $\tilde{\cal F}_{x_0+z}$\,-\,progressively measurable version of the process
$\{(U(x_0+z,\cdot)$, $\partial_1 U(x_0+z,\cdot));\,z\ge 0\}$
taking values in $E_1\times E_2$ is a stationary homogeneous strong Markov process
which is associated with the martingale problem of Definition \ref{our mart problem}
via a pathwise unique stochastic differential equation in $E_1\times E_2$
which can be formally written as
\begin{eqnarray*}
\dd u_z&=&v_z\,\dd z\\
\dd v_z&=&\rule{0pt}{20pt}
-\left[\rule{0pt}{11pt}\right.
(-\partial_t^2)^{\frac{1}{2}}u_z^a
+\sqrt{2}\,(-\partial_t^2)^{\frac{1}{4}}v_z^a
\left.\rule{0pt}{11pt}\right]\dd z
\,-\,\dd{\mathscr W}_z
\end{eqnarray*}
where $\{{\mathscr W}_z;\,z\ge 0\}$ stands for a ${\mathscr D}'$\,-\,valued Wiener process.
\end{theo}
\begin{cor}\label{corMtheo}
\begin{itemize}\item[(i)]
The unique weak solution $U(x,t)$ to (\ref{our SPDE}),
when seen as a process $U(x,\cdot)$ indexed by $x\ge x_0$ taking values in $E_1$,
satisfies
$$\pp\left[\left.\rule{0pt}{11pt}
\{U(\xi+y,\cdot)\in\Gamma\}\right|\tilde{\cal F}_\xi
\right]
\,\stackrel{\mbox{\tiny\rm a.s.}}{=}\,
\pp\left[\left.\rule{0pt}{11pt}
\{U(\xi+y,\cdot)\in\Gamma\}\right|\sigma(U(\xi,\cdot),\partial_1 U(\xi,\cdot))
\right]$$
for any finite $\tilde{\cal F}_x$\,-\,stopping time $\xi\ge x_0$
and any $y\ge 0,\,\Gamma\in{\mathscr B}(E_1)$.
This remains valid when the filtration 
\rule{0pt}{12pt}$\tilde{\cal F}_x,\,x\ge x_0$,
is replaced by the filtration generated by the process
$\{(U(x,\cdot),\partial_1 U(x,\cdot));\,x\ge x_0\}$
augmented by the $\pp$-null sets in ${\cal F}$.
\item[(ii)]
For any $x\in\NR$, the $\sigma$-algebras
$\sigma\{U(y,t):y<x,t>0\}$ and $\sigma\{U(y,t):y>x,t>0\}$
are conditionally independent given $\sigma(U(x,\cdot),\partial_1 U(x,\cdot))$
where
$$\sigma(U(x,\cdot),\partial_1 U(x,\cdot))
\;\varsubsetneq\;
germ\left(\rule{0pt}{10pt}\right.
\{x\}\times(0,\infty)
\left.\rule{0pt}{10pt}\right)
\;\;\stackrel{\mbox{\tiny def}}{=}\hspace{-10pt}
\bigcap_{
\substack{O\,open\,in\,{\cal Q}_+\\ 
\rule{0pt}{9pt}
\{x\}\times(0,\infty)\subseteq O}
}\hspace{-10pt}
\sigma\{U(y,t):(y,t)\in O\}.$$
\end{itemize}
\end{cor}
\section{Proofs}\label{proofs}
{\it Proof} of {\bf Proposition \ref{stateSpace}}.
For (i) fix $(x,h)\in[x_0,\infty)\times{\mathscr D}$ and notice that
$$\int_{x_0}^x\partial_1 U(y,h)\,\dd y\,=\,
\int_\NR(
\int\hspace{-10pt}\int_{{\cal Q}_+}\hspace{-10pt}B(\dd y',\dd s')
\left[\rule{0pt}{10pt}\right.\hspace{-3pt}
\int_0^\infty\hspace{-10pt}\partial_3 g(y',s'\,;y,t){h}(t)\,\dd t
\hspace{-3pt}\left.\rule{0pt}{10pt}\right])\,\ind_{(x_0,x]}(y)\,\dd y$$
by the definition of $\partial_1 U(y,h)$.
The integral on the right-hand side a.s.\ equals
\begin{equation}\label{int rhs}
\int\hspace{-10pt}\int_{{\cal Q}_+}\hspace{-10pt}B(\dd y',\dd s')
\int_\NR
\left[\rule{0pt}{10pt}\right.\hspace{-3pt}
\int_0^\infty\hspace{-10pt}\partial_3 g(y',s'\,;y,t){h}(t)\dd t
\hspace{-3pt}\left.\rule{0pt}{10pt}\right]\,\ind_{(x_0,x]}(y)\,\dd y
\end{equation}
by applying Lemma \ref{fubini} with respect to the bounded function $\phi=\ind_{(x_0,x]}$.
Here the condition of Lemma \ref{fubini}(ii) is easily satisfied
because the covariance of $\partial_1 U(y,h)$ given in Proposition \ref{stateSpace}(ii)
does not depend on $y$.
Then the equation for $U(x,h)$ follows from (\ref{int rhs})
by Fubini's theorem with respect to $\dd t\dd y$
which can be applied for every $(y',s')\in{\cal Q}_+$ because
$$\int_{x_0}^x
\int_0^\infty\hspace{-0pt}|\partial_3 g(y',s'\,;y,t){h}(t)|\,\dd t\dd y\,<\,\infty.$$

In order to show the equation for $\partial_1 U(x,h)$ we first calculate:
{\small
\begin{eqnarray*}
&&\partial_1 U(x,h)-\partial_1 U(x_0,h)\\
&=&\rule{0pt}{20pt}
\int\hspace{-10pt}\int_{{\cal Q}_+}\hspace{-10pt}B(\dd y',\dd s')
\left[\rule{0pt}{10pt}\right.\hspace{-3pt}
\int_0^\infty\hspace{-10pt}\partial_3 g(y',s'\,;x,t){h}(t)\,\dd t
\hspace{-3pt}\left.\rule{0pt}{10pt}\right]
-
\int\hspace{-10pt}\int_{{\cal Q}_+}\hspace{-10pt}B(\dd y',\dd s')
\left[\rule{0pt}{10pt}\right.\hspace{-3pt}
\int_0^\infty\hspace{-10pt}\partial_3 g(y',s'\,;x_0,t){h}(t)\,\dd t
\hspace{-3pt}\left.\rule{0pt}{10pt}\right]\\
&\stackrel{\mbox{\tiny\rm a.s.}}{=}&\rule{0pt}{20pt}
\int\hspace{-10pt}\int_{{\cal Q}_+}\hspace{-10pt}B(\dd y',\dd s')
\left[\rule{0pt}{10pt}\right.\hspace{-3pt}
\int_{s'}^\infty\hspace{-5pt}
[\int_{x_0}^x 
\partial_3^2 g(y',s'\,;y,t)\,\dd y]\,{h}(t)\,\dd t
\hspace{-3pt}\left.\rule{0pt}{10pt}\right]\\
&=&\rule{0pt}{20pt}
\int\hspace{-10pt}\int_{{\cal Q}_+}\hspace{-10pt}B(\dd y',\dd s')
\left[\rule{0pt}{10pt}\right.\hspace{-3pt}
\int_{s'}^\infty\hspace{-5pt}
[\int_{x_0}^x
\partial_4 g(y',s'\,;y,t)\,\dd y]\,{h}(t)\,\dd t
\hspace{-3pt}\left.\rule{0pt}{10pt}\right]\\
&=&\rule{0pt}{20pt}
\int\hspace{-10pt}\int_{{\cal Q}_+}\hspace{-10pt}B(\dd y',\dd s')
\left[\rule{0pt}{10pt}\right.\hspace{-3pt}
\int_{s'}^\infty\hspace{-5pt}
\partial_t[\int_{x_0}^x
g(y',s'\,;y,t)\,\dd y]\,{h}(t)\,\dd t
\hspace{-3pt}\left.\rule{0pt}{10pt}\right]\\
&\stackrel{\mbox{\tiny\rm a.s.}}{=}&\rule{0pt}{20pt}
\hspace{-10pt}
-\int\hspace{-10pt}\int_{{\cal Q}_+}\hspace{-10pt}B(\dd y',\dd s')
\left[\rule{0pt}{10pt}\right.\hspace{-3pt}
\int_{s'}^\infty\hspace{-5pt}
[\int_{x_0}^x
g(y',s';y,t)\,\dd y]\,{h}'(t)\,\dd t
\hspace{-3pt}\left.\rule{0pt}{10pt}\right]
-
\int\hspace{-10pt}\int_{{\cal Q}_+}\hspace{-10pt}B(\dd y',\dd s')
\lim_{t\downarrow s'}[\int_{x_0}^x g(y',s';y,t)\,\dd y]\,h(s').
\end{eqnarray*}
}\noindent
Again applying Fubini's Theorem and our stochastic Fubini Lemma \ref{fubini},
one sees that
$$\int\hspace{-10pt}\int_{{\cal Q}_+}\hspace{-10pt}B(\dd y',\dd s')
\left[\rule{0pt}{10pt}\right.\hspace{-3pt}
\int_{s'}^\infty\hspace{-5pt}
[\int_{x_0}^x
g(y',s'\,;y,t)\,\dd y]\,{h}'(t)\,\dd t
\hspace{-3pt}\left.\rule{0pt}{10pt}\right]
\,=\,
\int_{x_0}^x U(y,h')\,\dd y$$
hence the equation for $\partial_1 U(x,h)$ follows since
$$\int\hspace{-10pt}\int_{{\cal Q}_+}\hspace{-10pt}B(\dd y',\dd s')
\lim_{t\downarrow s'}[\int_{x_0}^x g(y',s'\,;y,t)\,\dd y]\,h(s')
\,\stackrel{\mbox{\tiny\rm a.s.}}{=}\,
\dB\,
(\ind_{(x_0,x]}\otimes h)(y,s)$$
holds true by the strong continuity of the heat semigroup in $L^2(\NR)$.

For proving item (ii) of the proposition
fix $x\ge x_0$ and $h_1,h_2\in{\mathscr D}$. Then
{\small
\begin{eqnarray*}
&&\ee\;U(x,h_1)U(x,h_2)\\
&=&\rule{0pt}{20pt}
\ee
\dB
\left[\rule{0pt}{10pt}\right.\hspace{-3pt}
\int_0^\infty\hspace{-10pt}\,g(y,s\,;x,t){h}_1(t)\,\dd t
\hspace{-3pt}\left.\rule{0pt}{10pt}\right]
\int\hspace{-10pt}\int_{{\cal Q}_+}\hspace{-10pt}B(\dd y',\dd s')
\left[\rule{0pt}{10pt}\right.\hspace{-3pt}
\int_0^\infty\hspace{-10pt}\,g(y',s'\,;x,t'){h}_2(t')\,\dd t'
\hspace{-3pt}\left.\rule{0pt}{10pt}\right]\\
&=&\rule{0pt}{20pt}
\int_\NR\int_0^\infty
\left[\rule{0pt}{10pt}\right.\hspace{-3pt}
\int_0^\infty\hspace{-10pt}\,g(y,s\,;x,t){h}_1(t)\,\dd t
\hspace{-3pt}\left.\rule{0pt}{10pt}\right]
\left[\rule{0pt}{10pt}\right.\hspace{-3pt}
\int_0^\infty\hspace{-10pt}\,g(y,s\,;x,t'){h}_2(t')\,\dd t'
\hspace{-3pt}\left.\rule{0pt}{10pt}\right]\,\dd s\dd y\\
&=&\rule{0pt}{20pt}
\int_0^\infty h_1(t)
\int_0^\infty
\left[\rule{0pt}{10pt}\right.\hspace{-3pt}
\int_\NR\int_0^\infty\hspace{-5pt}g(y,s\,;x,t)g(y,s\,;x,t')\,\dd s\dd y
\hspace{-2pt}\left.\rule{0pt}{10pt}\right]
h_2(t')\,\dd t'\,\dd t\\
&=&
\langle\,h_1\,;\frac{-\sqrt{|\cdot|}}{\sqrt{4\pi}}\ast h_2^a\,\rangle
\end{eqnarray*}
}\noindent
because
$$\frac{1}{\sqrt{4\pi}}\,(\sqrt{t+t'}-\sqrt{|t-t'|})
\,=\,
\int_\NR\int_0^\infty\hspace{-5pt}g(y,s\,;x,t)g(y,s\,;x,t')\,\dd s\dd y.$$
We only mention that, by the well-known properties of the Green's function $g$,
the integrability conditions needed for the above calculation are satisfied
in the case of test functions $h_1,h_2$ with \underline{compact} support.

The covariance of the process $\{\partial_1 U(x,h);\,h\in{\mathscr D}\}$
can be verified by a similar calculation since
$$\frac{1}{2\sqrt{4\pi}}\,
(\frac{1}{\sqrt{|t-t^{'}|}}-\frac{1}{\sqrt{t+t^{'}}})
\,=
\int_\NR\int_0^\infty\partial_3 g(y,s\,;x,t)  \partial_3 g(y,s\,;x,t')\,\dd s\dd y$$
and 
$$\int_\NR\int_0^\infty g(y,s\,;x,t)  \partial_3 g(y,s\,;x,t')\,\dd s\dd y
\,=\,0$$
for all $t,t'\ge 0$ gives the independence of the two processes.

We now show part (iii) of the proposition.
Fix $x\ge x_0$ and $\alpha>3/2$.
If $h\in C_{0,\alpha}$ then
{\small
\begin{eqnarray*}
&&\int_\NR\int_0^\infty\int_0^\infty g(y,s\,;x,t)^2\,|h(t)|\,\dd t\,\dd s\dd y\\
&\le&\rule{0pt}{20pt}
\int_\NR\int_0^\infty\int_s^\infty
\frac{1}{4\pi(t-s)}\exp\{\frac{-(x-y)^2}{2(t-s)}\}(1+t)^{-\alpha}\,\dd t\,\dd s\dd y
\cdot\|h\|_{0,\alpha}\\
&=&\rule{0pt}{20pt}
\frac{\|h\|_{0,\alpha}}{2\sqrt{2\pi}}
\int_0^\infty\int_s^\infty(t-s)^{-\frac{1}{2}}
(1+t)^{-\alpha}\,\dd t\dd s
\,=\,
\frac{\|h\|_{0,\alpha}}{\sqrt{2\pi}}
\int_0^\infty\sqrt{t}\,(1+t)^{-\alpha}\,\dd t\,<\,\infty
\end{eqnarray*}
}\noindent
because $\alpha>3/2$.
Hence, by Lemma \ref{fubini}, 
the process $\{U(x,h);\,h\in{\mathscr D}\}$
can be extended to $h\in C_{0,\alpha}$ and
\begin{equation}\label{equal to U}
U(x,h)\,=\int_0^\infty U(x,t)\,h(t)\,\dd t
\quad\mbox{a.s. for all $h\in C_{0,\alpha}$.}
\end{equation}
But the above calculation also shows that
\begin{equation}\label{square of U}
\ee\int_0^\infty U(x,t)^2\,\,(1+t)^{-\alpha}\dd t\,<\,\infty
\end{equation}
thus
$$\int_0^\infty|U(x,t)|\,\,(1+t)^{-\alpha}\dd t\,<\,\infty
\quad\mbox{a.s.}$$
which implies
\begin{equation}\label{norm E1}
|\hspace{-5pt}\int_0^\infty U(x,t)\,h(t)\,\dd t|
\,\le\,
\|h\|_{0,\alpha}
\int_0^\infty|U(x,t)|\,\,(1+t)^{-\alpha}\dd t,
\quad\forall\,h\in C_{0,\alpha},\;\mbox{a.s.}
\end{equation}
As a consequence, 
there is a version of the process $\{U(x,h);\,h\in{\mathscr D}\}$ 
taking values in $(C_{0,\alpha})'$.
Since $U$ given by (\ref{greenrep}) is continuous in $(x,t)\in{\cal Q}_+$ such that
$\lim_{t\downarrow 0}U(x,t)=0$ for all $x\ge x_0$,
this version takes values in $E_1$ even. 

Next we prove item (iv) of Proposition \ref{stateSpace}.
Note that the Sobolev space $H_{\beta}$ can be identified with
$(Id-\partial_t^2)^{-\beta/2}L^2(\NR)$ in the sense of generalised functions.
Fix $\beta>1/4$ and define the operator
$Kh=(Id-\partial_t^2)^{-\beta/2}(\mathfrak{w}^{-1}h^a),\,h\in{\mathscr D}$.
Of course, $K^{-1}$ exists and it holds that
$K^{-1}h=\mathfrak{w}[(Id-\partial_t^2)^{\beta/2}h^a],\,h\in{\mathscr D}$.
Hence, if $v$ is a linear form with domain of definition which contains
$\{K e_i\}_{i=1}^\infty$ where
$\{e_i\}_{i=1}^\infty\subseteq{\mathscr D}$
is an orthonormal basis of $L^2([0,\infty))$
then
$$\sum_{i=1}^\infty|v(K e_i)|^2\,<\,\infty
\quad\mbox{gives}\quad
v\,=\sum_{i=1}^\infty v(K e_i)
K^{-1}e_i\,\in\,(H_{\mathfrak{w},\beta}^a)'.$$

Fix $x\ge x_0$ and choose an orthonormal basis 
$\{e_i\}_{i=1}^\infty\subseteq{\mathscr D}$
of $L^2([0,\infty))$.
The above implies that if
the linear form $\partial_1 U(x,\cdot)$ can be extended to the linear hull of
${\mathscr D}\cup\{K e_i\}_{i=1}^\infty$
and if
\begin{equation}\label{conv in sobolev}
\ee\sum_{i=1}^\infty|\partial_1 U(x,K e_i)|^2\,<\,\infty
\end{equation}
then
\begin{equation}\label{version for partial_1 U}
\omega\,\mapsto
\ind_{\{\sum_{i=1}^\infty|\partial_1 U(x,K e_i)|^2<\infty\}}(\omega)
\sum_{i=1}^\infty\partial_1 U(\omega,x,K e_i)
K^{-1}e_i\end{equation}
defines a version of $\partial_1 U(x,\cdot)$ taking values in $(H_{\mathfrak{w},\beta}^a)'$.

In order to show (\ref{conv in sobolev}) recall from Remark \ref{why H minus two}(ii) that
$C_2 h=const(-\partial_t^2)^{-\frac{1}{4}}\,h^a$. First,
applying Proposition \ref{stateSpace}(ii), 
we have that\footnote{See Section \ref{notation} for  
$\|l^a\|^2_{L^2(\NR)}=2\langle l\,;l\rangle$.}
\begin{eqnarray*}
\ee\,|\partial_1 U(x,K e_i)|^2
&=&
\langle(Id-\partial_t^2)^{-\beta/2}(\mathfrak{w}^{-1}e_i^a)\,;
C_2(Id-\partial_t^2)^{-\beta/2}(\mathfrak{w}^{-1}e_i^a)\rangle\\
&=&\rule{0pt}{20pt}\frac{const}{2}\,
\|(-\partial_t^2)^{-\frac{1}{8}}(Id-\partial_t^2)^{-\beta/2}(\mathfrak{w}^{-1}e_i^a)\|^2_{L^2(\NR)}\\
&=&\rule{0pt}{20pt}\frac{const}{2\cdot 2\pi}\,
\left\|\rule{0pt}{11pt}\right.|\cdot|^{-1/4}(1+|\cdot|^2)^{-\beta/2}(\mathfrak{w}^{-1}e_i^a)^F
\left.\rule{0pt}{11pt}\right\|^2_{L^2(\NR)}\\
&\le&\rule{0pt}{20pt}\frac{const}{2\cdot 2\pi}\,
\left\|\rule{0pt}{11pt}\right.|\cdot|^{-1/4}(\mathfrak{w}^{-1}e_i^a)^F
\left.\rule{0pt}{11pt}\right\|^2_{L^2(\NR)}
\,\le\,
\langle\,|\mathfrak{w}^{-1}e_i|\,;C_2|\mathfrak{w}^{-1}e_i|\,\rangle\,<\,\infty
\end{eqnarray*}
for each single $i$ since $|\mathfrak{w}^{-1}e_i|\le|e_i|$ and
$e_i\in{\mathscr D}$. As a consequence,
$\partial_1 U(x,\cdot)$ can be extended to the linear hull of
${\mathscr D}\cup\{K e_i\}_{i=1}^\infty$
and the left-hand side of (\ref{conv in sobolev}) makes sense.

Taking into account the calculations of the last paragraph, 
condition (\ref{conv in sobolev}) becomes equivalent to
$$\sum_{i=1}^\infty
\left\|\rule{0pt}{11pt}\right.
|\cdot|^{-1/4}(1+|\cdot|^2)^{-\beta/2}(\mathfrak{w}^{-1}e_i^a)^F
\left.\rule{0pt}{11pt}\right\|^2_{L^2(\NR)}
\,<\,\infty$$
where
$$(\mathfrak{w}^{-1}e_i^a)^F\,=\,
{\bf i}\int_\NR\sin(-\tau t)\,\mathfrak{w}^{-1}(t)e_i^a(t)\,\dd t
\,=\,
-2{\bf i}\int_0^\infty\sin(\tau t)\,\mathfrak{w}^{-1}(t)e_i(t)\,\dd t.$$
Thus
\begin{eqnarray*}
&&\sum_{i=1}^\infty
\left\|\rule{0pt}{11pt}\right.
|\cdot|^{-1/4}(1+|\cdot|^2)^{-\beta/2}(\mathfrak{w}^{-1}e_i^a)^F
\left.\rule{0pt}{11pt}\right\|^2_{L^2(\NR)}\\
&=&\rule{0pt}{20pt}
\sum_{i=1}^\infty\int_\NR
\frac{1}{\sqrt{|\tau|}\,(1+|\tau|^2)^{\beta}}
\int_0^\infty\sin(\tau t)\,\mathfrak{w}^{-1}(t)e_i(t)\,\dd t
\int_0^\infty\sin(\tau t')\,\mathfrak{w}^{-1}(t')e_i(t')\,\dd t'
\,\dd\tau\\
&=&\rule{0pt}{30pt}
\int_0^\infty\sum_{i=1}^\infty
\left\{\rule{0pt}{11pt}\right.
\int_0^\infty
[\int_\NR
\frac{\sin(\tau t)\,\mathfrak{w}^{-1}(t)\sin(\tau t')\,\mathfrak{w}^{-1}(t')}
{\sqrt{|\tau|}\,(1+|\tau|^2)^{\beta}}
\,\dd\tau\,]
\,e_i(t')\,\dd t'
\left.\rule{0pt}{11pt}\right\}
e_i(t)\,\dd t\\
&=&\rule{0pt}{30pt}
\int_0^\infty\!\!\int_\NR
\frac{\sin(\tau t)^2\,\mathfrak{w}^{-1}(t)^2}
{\sqrt{|\tau|}\,(1+|\tau|^2)^{\beta}}
\,\dd\tau\dd t
\,\le\,
\int_\NR
\frac{\dd\tau}
{\sqrt{|\tau|}\,(1+|\tau|^2)^{\beta}}
\;
\int_0^\infty\!\!
\frac{\dd t}
{(1+|t|^{\frac{1}{2}+\varepsilon})^2}
\,<\,\infty
\end{eqnarray*}
since $\beta>1/4$.

We finally justify item (v) of Proposition \ref{stateSpace}.
First, the stationarity follows from item (ii) because the covariances
do not depend on $x\ge x_0$. Second, $(\omega,x,t)\mapsto U(\omega,x,t)$ is clearly 
${\cal F}\otimes{\mathscr B}([x_0,\infty))\otimes{\mathscr B}([0,\infty))$\,-\,measurable
leading to an ${\cal F}\otimes{\mathscr B}([x_0,\infty))$\,-\,measurable version
of $x\mapsto U(x,\cdot)\in E_1$
and the version of $(\omega,x)\mapsto\partial_1 U(\omega,x,\cdot)$
given by (\ref{version for partial_1 U}) is also
${\cal F}\otimes{\mathscr B}([x_0,\infty))$\,-\,measurable.
Third, using stationarity, (\ref{bound on norms}) already follows from
$$\ee\|U(x,\cdot)\|_{E_1}^2\,<\,\infty
\quad\mbox{and}\quad
\ee\|\partial_1 U(x,\cdot)\|_{E_2}^2\,<\,\infty$$
for an arbitrary but fixed $x\ge x_0$
where the first expectation is finite because of (\ref{square of U}),(\ref{norm E1})
and the second expectation is equal to the left-hand side of (\ref{conv in sobolev})
which was shown to be finite above.
Finally, the existence of a continuous version on the closure of ${\cal Q}_+^0$
of the solution $(x,t)\mapsto U(x,t)$ as given by (\ref{greenrep}) 
is standard--see \cite{W1986}.\hfill$\Box$\\

{\it Proof of} {\bf Proposition \ref{new wipro}}.
Recalling Remark \ref{idea new wipro}(i), we only deal with the following two issues
and refer to Theorem 12.1 in \cite{Y1997} otherwise.

First, after several steps, one has to identify
$$\int_{x_0}^{x_0+\,\fatdot}
\ee\left[\left.\rule{0pt}{11pt}
F(L)\varrho_l(L,y)\right|{\cal F}_{(-\infty,y]}
\right]\dd y$$
with the covariation between the martingales
$$\ee\left[\left.\rule{0pt}{11pt}
F(L)\right|{\cal F}_{(-\infty,x_0+\,\fatdot\,]}
\right]
\quad\mbox{and}\quad
\dB\,(\ind_{(x_0,x_0+\,\fatdot\,]}\otimes l)(y,s).$$
But, by the assumptions on $\varrho_l$ made in the proposition, it holds that
$$\int_{x_0}^{x_0+z}
\ee\left[\left.\rule{0pt}{11pt}
F(L)\varrho_l(L,y)\right|{\cal F}_{(-\infty,y]}
\right]\dd y
\,=\,
\int_{x_0}^{x_0+z}
\int_0^\infty\dot{\lambda}_{ys}(F)\,l(s)\,\dd s\,\dd y,
\quad\mbox{$z\ge 0$, a.s.,}$$
where, by the martingale representation of $F(L)$,
the above right-hand side is the wanted covariation.

Second, 
if $\varrho_l$ exists for $l_1,l_2\in L^2([0,\infty))$ and if $a_1,a_2\in\NR$ then
$$\int_0^\infty\dot{\lambda}_{ys}(F)\,(a_1 l_1+a_2 l_2)(s)\,\dd s
\,\stackrel{\mbox{\tiny\rm a.s.}}{=}\,
\ee\left[\left.\rule{0pt}{11pt}
F(L)(a_1\varrho_{l_1}(L,y)+a_2\varrho_{l_2}(L,y))\right|{\cal F}_{(-\infty,y]}
\right].$$
But 
$(\omega,\phi,y)\mapsto 
a_1\varrho_{l_1}(\omega,\phi,y)+a_2\varrho_{l_2}(\omega,\phi,y)$
also satisfies the other properties of $\varrho_l$ stated in the proposition
hence it can be taken to be $\varrho_{a_1 l_1+a_2 l_2}$.
Then the linearity (\ref{lin}) follows from the uniqueness of $\varrho_l$ 
explained in Remark \ref{idea new wipro}(ii).
Note that (\ref{lin})
is not required for the argument given in Remark \ref{idea new wipro}(ii).
\hfill$\Box$\\

{\it Proof of} {\bf Lemma \ref{identify rho}}.
Fix $\nu>0$ and observe that, by change of variable $(t'=tr)$,
the test function $l_\nu$ can be represented as
\begin{equation}\label{rep l_nu}
l_\nu(t)\,=\,\frac{\sqrt{t}}{\sqrt{4\pi}}\int_0^\infty
(\frac{1}{\sqrt{|1-r|}}-\frac{1}{\sqrt{1+r}})\,e^{-\nu tr}\,\dd r,\quad t\ge 0.
\end{equation}
Thus, since the function 
$r\mapsto\left|\rule{0pt}{10pt}\right.
|1-r|^{-1/2}-(1+r)^{-1/2}
\left.\rule{0pt}{10pt}\right|^p$ on $[0,\infty)$
is integrable for all $1\le p<2$, 
Lebesgue's dominated convergence theorem implies both the continuity of $l_\nu$
and $l_\nu(t)\to 0$ if $t$ tends to zero.
Moreover, $l_\nu\in L^2([0,\infty))\cap L^\infty([0,\infty))$ follows from
\begin{eqnarray*}
&&\int_0^\infty
(\frac{1}{\sqrt{|1-r|}}-\frac{1}{\sqrt{1+r}})\,e^{-\nu tr}\,\dd r\\
&\le&\rule{0pt}{20pt}
e^{-\nu t}\int_0^1 
(\frac{1}{\sqrt{r}}-\frac{1}{\sqrt{2-r}})\,e^{\nu tr}\,\dd r
\,+\,
\int_0^\infty
\frac{1}{\sqrt{r}}\,e^{-\nu t(r+1)}\,\dd r\\
&\le&\rule{0pt}{20pt}
e^{-\nu t}
\left[\rule{0pt}{11pt}\right.\!
\sqrt{2}\,e^{\nu t/2}
+
\hspace{-2.5cm}\underbrace{
\int_{1/2}^1 
(\frac{1}{\sqrt{r}}-\frac{1}{\sqrt{2-r}})\,e^{\nu tr}\,\dd r
}_{\mbox{$
-(\sqrt{2}-\sqrt{\frac{2}{3}}\,)\,\frac{e^{\nu t/2}}{\nu t}
+
\frac{1}{2}\int_{1/2}^1\,
(r^{-3/2}+(2-r)^{-3/2})\,\frac{e^{\nu tr}}{\nu t}\,\dd r
$}}\hspace{-2.3cm}
\!\left.\rule{0pt}{11pt}\right]
+e^{-\nu t}\sqrt{\pi/(\nu t)}
\end{eqnarray*}
since
$$\int_{1/2}^1\,
(r^{-3/2}+(2-r)^{-3/2})\,\frac{e^{\nu tr}}{\nu t}\,\dd r
\,\le\,
(2^{3/2}+1)\,[\,e^{\nu t}-e^{\nu t/2}\,]\,
\frac{1}{(\nu t)^2}$$
which, in the end, yields
$l_\nu(t)={\cal O}(t^{-3/2})$ for $t\to\infty$ by (\ref{rep l_nu}).

The next step is to identify $\varrho_{l_\nu}(\phi,y)$ as given in the lemma
so, in particular, we have to show (\ref{how rho}). 
Fix $y>x_0$ and $F\in{\mathfrak F}C_b^\infty({\mathscr D})$ given by
$$F(\phi)\,=\,f(\langle\phi\,;h_1\rangle,\dots,\langle\phi\,;h_n\rangle),\quad\phi\in E_1,$$
where $f\in C_b^\infty(\NR^n),\,h_i\in{\mathscr D},\,i=1,\dots,n$, for some $n\ge 1$.

Since $U(x_0,\cdot)$ is a stochastic integral against the Brownian sheet $B$
with deterministic integrand,
the Malliavin derivative $D_{ys}F(U(x_0,\cdot))$ exists and can explicitly be given by
$$D_{ys}F(U(x_0,\cdot))\,=\sum_{i=1}^n
\partial_i f(\dots,\langle U(x_0,\cdot)\,;h_i\rangle,\dots)
\int_0^\infty g(y,s\,;x_0,t)h_i(t)\,\dd t.$$
Furthermore, by Clark-Ocone's formula, we have the identity
$$\dot{\lambda}_{ys}(F)\,=\,
\ee\left[\left.\rule{0pt}{11pt}
D_{ys}F(U(x_0,\cdot))\right|{\cal F}_{(-\infty,y]}
\right],
\quad\mbox{$s\ge 0$, a.s.}$$
Therefore we obtain that
\begin{eqnarray*}
&&\int_0^\infty\dot{\lambda}_{ys}(F)\,l_\nu(s)\,\dd s\\
&\stackrel{\mbox{\tiny\rm a.s.}}{=}&
\sum_{i=1}^n\ee\left[\left.\rule{0pt}{11pt}
\partial_i f(\dots,\langle U(x_0,\cdot)\,;h_i\rangle,\dots)\right|{\cal F}_{(-\infty,y]}
\right]
\int_0^\infty[\int_0^\infty g(y,s\,;x_0,t)l_\nu(s)\,\dd s\,]\,h_i(t)\,\dd t\\
&=&
\sum_{i=1}^n\ee\left[\left.\rule{0pt}{11pt}
\partial_i f(\dots,\langle U(x_0,\cdot)\,;h_i\rangle,\dots)\right|{\cal F}_{(-\infty,y]}
\right]
\langle g_y^{x_0}\ast l_\nu^0\,;h_i\rangle
\end{eqnarray*}
where 
$l_\nu^0\stackrel{\mbox{\tiny def}}{=}l_\nu\ind_{(0,\infty)}$
is treated as a function on $\NR$
and $g_y^{x_0}$ was defined in Section \ref{notation}.

Since $U(x_0,\cdot)-U(x_0,\cdot)_y$ and ${\cal F}_{(-\infty,y]}$ are independent,
the last sum of conditional expectations simplifies to
$$\int_{E_1}\sum_{i=1}^n
\partial_i f(\dots,\langle U(x_0,\cdot)_y\,;h_i\rangle+\langle\phi^y\,;h_i\rangle,\dots)
\langle g_y^{x_0}\ast l_\nu^0\,;h_i\rangle\,\mu^y(\dd\phi^y)$$
where $\mu^y$ denotes the image measure of $U(x_0,\cdot)-U(x_0,\cdot)_y$
on $E_1$ equipped with the Borel-$\sigma$-algebra.
Remark that, similar to the proof of Proposition \ref{stateSpace}(iii),
there are versions of $U(x_0,\cdot)_y$ and $U(x_0,\cdot)-U(x_0,\cdot)_y$
taking values in $E_1$.

Now introduce the function
$$F_y(\omega,\phi^y)\,\stackrel{\mbox{\tiny def}}{=}\,
F(U(x_0,\cdot)_y(\omega)+\phi^y)$$
and observe that
$$\sum_{i=1}^n
\partial_i f(\dots,\langle U(x_0,\cdot)_y\,;h_i\rangle+\langle\phi^y\,;h_i\rangle,\dots)
\langle g_y^{x_0}\ast l_\nu^0\,;h_i\rangle
\,=\,
\frac{\partial F_y(\phi^y)}{\partial(g_y^{x_0}\ast l_\nu^0)}$$
where the right-hand side is the G\^{a}teaux derivative
into the direction $g_y^{x_0}\ast l_\nu^0$ defined by
$$\frac{\dd}{\dd r}\,F_y\!\left.\left(\rule{0pt}{10pt}
\phi^y+r(g_y^{x_0}\ast l_\nu^0)\right)\right|_{r=0}.$$
Note that this requires $g_y^{x_0}\ast l_\nu^0\in E_1$ which can easily be verified
using the explicit structure of $g_y^{x_0}$ 
because $l_\nu^0$ is bounded.

Having found that
\begin{equation}\label{equal gateaux deriv}
\int_0^\infty\dot{\lambda}_{ys}(F)\,l_\nu(s)\,\dd s\\
\,\stackrel{\mbox{\tiny\rm a.s.}}{=}\,
\int_{E_1}\frac{\partial F_y(\phi^y)}{\partial(g_y^{x_0}\ast l_\nu^0)}\,\mu^y(\dd\phi^y)
\end{equation}
we now want to justify that
$$\int_{E_1}\frac{\partial F_y(\phi^y)}{\partial(g_y^{x_0}\ast l_\nu^0)}\,\mu^y(\dd\phi^y)
\,\stackrel{\mbox{\tiny\rm a.s.}}{=}\,
\ee\left[\left.\rule{0pt}{11pt}
F(U(x_0,\cdot))
\langle U(x_0,\cdot)-U(x_0,\cdot)_y\,;
\frac{2\sqrt{\nu}e^{-\nu\,\fatdot}}{e^{-\sqrt{\nu}(y-x_0)}}\rangle
\right|{\cal F}_{(-\infty,y]}
\right].$$
But
\begin{equation}\label{partial int}
\int_{E_1}\frac{\partial F_y(\phi^y)}{\partial(g_y^{x_0}\ast l_\nu^0)}\,\mu^y(\dd\phi^y)
\,=\,
\int_{E_1}F_y(\phi^y)\langle\phi^y\,;C_y^{-1}(g_y^{x_0}\ast l_\nu^0)\rangle\,\mu^y(\dd\phi^y)
\end{equation}
if the direction $g_y^{x_0}\ast l_\nu^0$ 
is in the Cameron-Martin space $H_y$ of the Gaussian measure $\mu^y$
with covariance $C_y:E_1'\to E_1''$.
\begin{rem}\label{details of method}\rm
\begin{itemize}\item[(i)]
We refer to \cite{B1998} being a good reference for the theory of Gaussian measures
on infinite-dimensional spaces. The covariance $C_y:E_1'\to E_1''$ can be extended
to the reproducing kernel Hilbert space $H_y'$ of $\mu^y$
and $C_y$ acts on $H_y'$ as an isomorphism between $H_y'$ and the
Cameron-Martin space $H_y,\,H_y\subseteq E_1\subseteq E_1''$.
So, checking if $g_y^{x_0}\ast l_\nu^0\in H_y$ can be done by finding 
a solution $m_\nu^y\in H_y'$ of the equation
\begin{equation}\label{crucial equation}
g_y^{x_0}\ast l_\nu^0\,=\,C_y\,m_\nu^y
\end{equation}
and this will be the next step of the proof.
\item[(ii)]
It is clear that we have chosen $l_\nu$ in a way that (\ref{crucial equation})
can be solved in $H_y'$. When applying our method with respect to other
linear SPDEs with additive Gaussian noise then one has to study an equation
of similar type, i.e.
$$g_y^{x_0}\ast l\,=\,C_y\,m_l^y,$$
where $g_y^{x_0}$ comes from the Green's function associated with the SPDE
and $C_y$ is the covariance of some Gaussian measure. The task is then
to identify the `good' test functions $l$ for which such an equation
can be solved.
\end{itemize}
\end{rem}

We will show that
$$g_y^{x_0}\ast l_\nu^0=
C_y\,e^{-\nu\,\fatdot}/\,\widehat{g_y^{x_0}}(\nu)$$
which also implies that
the direction $g_y^{x_0}\ast l_\nu^0$ must be in $H_y$ because 
$e^{-\nu\,\fatdot}\in E_1'\subseteq H_y'$
and $C_y:H_y'\to H_y$ is an isomorphism.

Let's show the claimed equality.
As $C_y$ is the covariance of the image measure of the random variable
$U(x_0,\cdot)-U(x_0,\cdot)_y$ taking values in $E_1$ given by
$$U(x_0,t)-U(x_0,t)_y\,=
\int\hspace{-10pt}\int_{{\cal Q}_+^{y}}\hspace{-0pt}B(\dd y',\dd s')\,
g(y',s'\,;x_0,t),\quad t\ge 0,$$
we obtain that
\begin{eqnarray*}
C_y\,e^{-\nu\,\fatdot}(t)&=&
\int_0^\infty
\left[\rule{0pt}{10pt}\right.\hspace{-3pt}
\int_y^\infty\hspace{-8pt}\int_0^\infty\hspace{-10pt}
\,g(y',s'\,;x_0,t)g(y',s'\,;x_0,t')\,\dd s'\dd y'
\hspace{-0pt}\left.\rule{0pt}{10pt}\right]
e^{-\nu t'}\,\dd t'\\
&=&
\int_y^\infty\hspace{-4pt}
\left(\rule{0pt}{10pt}\right.
g_{y'}^{x_0}\ast\{
\int_0^\infty g(y',\cdot\,;x_0,t')e^{-\nu t'}\,\dd t'\,\ind_{(0,\infty)}(\cdot)
\}\left.\rule{0pt}{10pt}\right)
(t)\,\dd y'
\end{eqnarray*}
for all $t\ge 0$. Thus
\begin{eqnarray*}
\widehat{C_y\,e^{-\nu\,\fatdot}}(\tilde{\nu})
&=&
\int_y^\infty\hspace{-4pt}
\widehat{g_{y'}^{x_0}}(\tilde{\nu})
\underbrace{
\int_0^\infty
\left(\rule{0pt}{10pt}\right.
g_{y'}^{x_0}\ast[e^{-\tilde{\nu}\,\fatdot}\,\ind_{(0,\infty)}(\cdot)]
\left.\rule{0pt}{10pt}\right)(t')
\,e^{-\nu t'}\,\dd t'
}_{\widehat{g_{y'}^{x_0}}(\nu)(\tilde{\nu}+\nu)^{-1}}
\,\dd y'\\
&=&
\int_y^\infty\hspace{-0pt}
\frac{e^{-(y'-x_0)(\sqrt{\tilde{\nu}}+\sqrt{\nu})}}{4\sqrt{\tilde{\nu}\nu}(\tilde{\nu}+\nu)}
\,\dd y'
\,=\,
\widehat{g_{y}^{x_0}}(\tilde{\nu})\widehat{g_{y}^{x_0}}({\nu})
\,\frac{1}{(\sqrt{\tilde{\nu}}+\sqrt{\nu})(\tilde{\nu}+\nu)}
\end{eqnarray*}
such that
$$\widehat{C_y\,e^{-\nu\,\fatdot}}(\tilde{\nu})/\widehat{g_{y}^{x_0}}({\nu})
\,=\,
\widehat{g_{y}^{x_0}}(\tilde{\nu})
\,\frac{1}{(\sqrt{\tilde{\nu}}+\sqrt{\nu})(\tilde{\nu}+\nu)}
\,=\,
\widehat{g_y^{x_0}\ast l_\nu^0}(\tilde{\nu})$$
for all $\tilde{\nu}>0$ proving
$$C_y\,e^{-\nu\,\fatdot}/\,\widehat{g_y^{x_0}}(\nu)
\,=\,g_y^{x_0}\ast l_\nu^0$$
in the end.

The above allows us to use
$$e^{-\nu\,\fatdot}/\,\widehat{g_y^{x_0}}(\nu)
\,=\,
\frac{2\sqrt{\nu}e^{-\nu\,\fatdot}}{e^{-\sqrt{\nu}(y-x_0)}}
\quad\mbox{for}\quad
C_y^{-1}(g_y^{x_0}\ast l_\nu^0)$$
on the right-hand side of (\ref{partial int}) leading to
\begin{eqnarray*}
\int_{E_1}\frac{\partial F_y(\phi^y)}{\partial(g_y^{x_0}\ast l_\nu^0)}\,\mu^y(\dd\phi^y)
&=&
\int_{E_1}F_y(\phi^y)\langle\phi^y\,;
\frac{2\sqrt{\nu}e^{-\nu\,\fatdot}}{e^{-\sqrt{\nu}(y-x_0)}}\rangle\,\mu^y(\dd\phi^y)\\
&\stackrel{\mbox{\tiny\rm a.s.}}{=}&\rule{0pt}{25pt}
\ee\left[\left.\rule{0pt}{11pt}
F(U(x_0,\cdot))
\langle U(x_0,\cdot)-U(x_0,\cdot)_y\,;
\frac{2\sqrt{\nu}e^{-\nu\,\fatdot}}{e^{-\sqrt{\nu}(y-x_0)}}\rangle
\right|{\cal F}_{(-\infty,y]}
\right].
\end{eqnarray*}
Because of (\ref{equal gateaux deriv}), this justifies that 
$\varrho_{l_\nu}(\phi,y)$ as given in Lemma \ref{identify rho} satisfies (\ref{how rho}).
It also satisfies the measurability conditions stated in Proposition \ref{new wipro}
and it only remains to show that 
$\varrho_{l_\nu}(U(x_0,\cdot),y)\in L^1(\Omega)$ for almost every $y\ge x_0$
and that
$y\mapsto\varrho_{l_\nu}(U(x_0,\cdot),y)$ is in $L^1([x_0,x])$ almost surely
for each $x\ge x_0$. But this follows from
\begin{eqnarray*}
&&\ee\int_{x_0}^x|\varrho_{l_\nu}(U(x_0,\cdot),y)|\,\dd y\\
&=&\ee\int_{x_0}^x|
\int\hspace{-10pt}\int_{{\cal Q}_+^{y}}\hspace{-0pt}B(\dd y',\dd s')
\,[
\int_0^\infty g(y',s'\,;x_0,t)\,
\frac{2\sqrt{\nu}e^{-\nu t}}{e^{-\sqrt{\nu}(y-x_0)}}\,\dd t
\,]\,|\,\dd y\\
&\le&\int_{x_0}^x\sqrt{
\int_\NR^\infty\hspace{-10pt}\int_0^\infty\hspace{-0pt}
[\int_0^\infty g(y',s'\,;x_0,t)\,
\frac{2\sqrt{\nu}e^{-\nu t}}{e^{-\sqrt{\nu}(y-x_0)}}\,\dd t\,]^2
\;\dd s'\dd y'}\;\dd y\\
&=&2\sqrt{\nu}\int_{x_0}^x\hspace{-3pt}e^{\sqrt{\nu}(y-x_0)}\,\dd y\hspace{3pt}
\sqrt{
\langle e^{-\nu\,\fatdot}\,;
\frac{-\sqrt{|\cdot|}}{\sqrt{4\pi}}\ast(e^{-\nu\,\fatdot})^a\rangle
}
\;<\,\infty
\end{eqnarray*}
where the equality in the last line is obtained 
by manipulations similar to the lines of proof of Proposition \ref{stateSpace}(ii).

We finally proof (\ref{as drift}). On the one hand, for fixed $y\ge x_0$, we have that
\begin{eqnarray*}
\varrho_{l_\nu}(U(x_0,\cdot),y)
&\stackrel{\mbox{\tiny\rm a.s.}}{=}&
\int\hspace{-10pt}\int_{{\cal Q}_+^{y}}\hspace{-0pt}B(\dd y',\dd s')
\int_0^\infty g(y',s'\,;x_0,t)\,
\frac{2\sqrt{\nu}e^{-\nu t}}{e^{-\sqrt{\nu}(y-x_0)}}\,\dd t\\
&=&
\int\hspace{-10pt}\int_{{\cal Q}_+^{y}}\hspace{-0pt}B(\dd y',\dd s')\,
2\sqrt{\nu}e^{\sqrt{\nu}(y-x_0)}
\underbrace{\int_{s'}^\infty g_{y'}^{x_0}(t-s')\,e^{-\nu t}\,\dd t}_{
e^{-\nu s'}\widehat{g_{y'}^{x_0}}(\nu)}\\
&=&
\int\hspace{-10pt}\int_{{\cal Q}_+^{y}}\hspace{-0pt}B(\dd y',\dd s')\,
e^{-\nu s'}e^{-\sqrt{\nu}(y'-y)}
\end{eqnarray*}
using again that $\widehat{g_{y'}^{x_0}}(\nu)=e^{-\sqrt{\nu}|y'-x_0|}/(2\sqrt{\nu})$
for $y'\in\NR$. On the other hand, it also holds that
\begin{eqnarray*}
&&U(y,\sqrt{\nu}e^{-\nu\,\fatdot})\,+\,\partial_1 U(y,e^{-\nu\,\fatdot})\\
&=&\rule{0pt}{20pt}\sqrt{\nu}
\int\hspace{-10pt}\int_{{\cal Q}_+}\hspace{-0pt}B(\dd y',\dd s')
\int_{s'}^\infty g(y',s'\,;y,t)\,e^{-\nu t}\,\dd t\\
&+&
\int\hspace{-10pt}\int_{{\cal Q}_+}\hspace{-0pt}B(\dd y',\dd s')
\int_{s'}^\infty\partial_3 g(y',s'\,;y,t)\,e^{-\nu t}\,\dd t
\end{eqnarray*}
\begin{eqnarray*}
&=&\sqrt{\nu}
\int\hspace{-10pt}\int_{{\cal Q}_+}\hspace{-0pt}B(\dd y',\dd s')\,
e^{-\nu s'}e^{-\sqrt{\nu}|y'-y|}/(2\sqrt{\nu})\\
&+&
\int\hspace{-10pt}\int_{{\cal Q}_+}\hspace{-0pt}B(\dd y',\dd s')
\underbrace{\int_{s'}^\infty
\frac{-2(y-y')}{4(t-s')\sqrt{4\pi(t-s')}}\,\exp\{\frac{-(y-y')^2}{4(t-s')}\}\,
e^{-\nu t}\,\dd t}_{
=\,\begin{array}[t]{l}
-\frac{1}{2}e^{-\nu s'}e^{-\sqrt{\nu}|y'-y|}\ind_{(-\infty,y]}(y')\\
+\frac{1}{2}e^{-\nu s'}e^{-\sqrt{\nu}|y'-y|}\ind_{(y,\infty)}(y')\rule{0pt}{15pt}
\end{array}}
\end{eqnarray*}
$$\stackrel{\mbox{\tiny\rm a.s.}}{=}\;
\int\hspace{-10pt}\int_{{\cal Q}_+^{y}}\hspace{-0pt}B(\dd y',\dd s')\,
e^{-\nu s'}e^{-\sqrt{\nu}(y'-y)}.\hspace{5.8cm}\eqno\Box$$
\begin{rem}\rm\label{partial_1 U adaptness}
The last part of the above proof also shows that
$$\partial_1 U(y,e^{-\nu\,\fatdot})\,\stackrel{\mbox{\tiny\rm a.s.}}{=}\,
U(y,\sqrt{\nu}e^{-\nu\,\fatdot})\,-
\int\hspace{-10pt}\int_{{\cal Q}_+\setminus{\cal Q}_+^{y}}\hspace{-0pt}B(\dd y',\dd s')\,
e^{-\nu s'}e^{-\sqrt{\nu}|y'-y|}$$
for all $y\ge x_0$.
\end{rem}
{\it Proof of} {\bf Proposition \ref{new drift}}.
Fix $h\in{\mathscr D}$.
Then $h^a$ is infinitely often differentiable but with compact support in $\NR$.
So, the function $({\mathfrak A}_2 h)^a$ defined by convolution is a $C^\infty$\,-\,function.

Next, choose an upper bound $c_h$ for the support of $h$ and fix $t>c_h$. Then
\begin{eqnarray*}
|{\mathfrak A}_2\,h(t)|&=&|\hspace{-3pt}
\int_{-c_h}^{c_h}\frac{(h^{a})'(t')\,\dd t'}{\sqrt{\pi(t-t')}}|
\,=\,\frac{1}{2}\,|\hspace{-3pt}
\int_{-c_h}^{c_h}\frac{(h^{a})(t')\,\dd t'}{\sqrt{\pi}(t-t')^{3/2}}\,|\\
&\le&\rule{0pt}{20pt}
\sup_{t'\ge 0}|h(t')|\,c_h(t-c_h)^{-3/2}
\,=\,{\cal O}(t^{-3/2}),\quad t\to\infty,
\end{eqnarray*}
which yields
${\mathfrak A}_2\,h\in C_{0,\alpha}$ for $0\le\alpha<3/2$.
Since $\partial_t^k{\mathfrak A}_2\,h=[sgn(\cdot)(\sqrt{\pi}|\cdot|)^{-\frac{1}{2}}]\ast\partial_t^{k+1}(h^a)$,
the claim that
$\partial_t^k{\mathfrak A}_2\,h\in C_{0,\alpha}$ for $0\le\alpha<k+3/2,\,k=1,2,\dots$,
can be shown exactly the same way only using
$\int_{-c_h}^{c_h}(t-t')^{-k-3/2}\,\dd t'\,\le\,2c_h(t-c_h)^{-k-3/2}$ instead.

For ${\mathfrak A}_2\,h\in H_{\mathfrak{w},\beta}^a$ recall that
$\mathfrak{w}$ is a smooth weight function such that, for some $\varepsilon>0$,
$\mathfrak{w}\ge 1+|\cdot|^{\frac{1}{2}+\varepsilon}$ 
but
$\mathfrak{w}=1+|\cdot|^{\frac{1}{2}+\varepsilon}$ outside a neighbourhood of zero.
Hence,
$\partial_t^k{\mathfrak A}_2\,h\in C_{0,\alpha}$ for $0\le\alpha<k+3/2$ implies
$\partial_t^k[\mathfrak{w}({\mathfrak A}_2\,h)^a]\in C^\infty(\NR)\cap L^2(\NR)$
if $\varepsilon<k+1/2$ for $k=0,1,2,\dots$
Note that ${\mathfrak A}_2\,h\in H_{\mathfrak{w},0}^a$ if and only if
$\mathfrak{w}({\mathfrak A}_2\,h)^a\in L^2(\NR)$
hence, assuming $\varepsilon<1/2$,
it remains to discuss the case $\beta>0$.
Denote by $\lceil\beta\rceil$ the smallest integer larger than $\beta$. Then
\begin{eqnarray*}
&&\|(1+|\cdot|)^{\beta/2}\,(\mathfrak{w}({\mathfrak A}_2\,h)^a)^F\|^2_{L^2(\NR)}\\
&=&\rule{0pt}{20pt}
\int_\NR(1+\tau^2)^\beta\,|(\mathfrak{w}({\mathfrak A}_2\,h)^a)^F(\tau)|^2\,\dd\tau
\,\le\,
2^{\lceil\beta\rceil-1}
\int_\NR(1+\tau^{2\lceil\beta\rceil})\,|(\mathfrak{w}({\mathfrak A}_2\,h)^a)^F(\tau)|^2\,\dd\tau\\
&=&\rule{0pt}{20pt}
2^{\lceil\beta\rceil-1}2\pi
\left\|\rule{0pt}{10pt}\right.
\mathfrak{w}({\mathfrak A}_2\,h)^a
\left.\rule{0pt}{10pt}\right\|^2_{L^2(\NR)}
+
2^{\lceil\beta\rceil-1}2\pi
\left\|\rule{0pt}{10pt}\right.
\partial_t^{\lceil\beta\rceil}[\mathfrak{w}({\mathfrak A}_2\,h)^a]
\left.\rule{0pt}{10pt}\right\|^2_{L^2(\NR)}
\,<\,\infty
\end{eqnarray*}
proving ${\mathfrak A}_2\,h\in H_{\mathfrak{w},\beta}^a$.
Note that $\varepsilon<1/2$ is needed
for the finiteness of the first summand in the last line, again.

Using the large\,-$t$-behaviour of $({\mathfrak A}_2\,h)'$ found above,
${\mathfrak A}_1{\mathfrak A}_2 h$ is well-defined and, 
for $\frac{1}{2}<\alpha'<\frac{5}{2}$,
we obtain that
\begin{eqnarray}
|{\mathfrak A}_1{\mathfrak A}_2 h(t)|
&\le&
\int_t^\infty\frac{|({\mathfrak A}_2 h)'(t')|\,\dd t'}{\sqrt{t'-t}}
\,\le\,
\|({\mathfrak A}_2\,h)'\|_{0,\alpha'}\int_0^\infty\frac{\dd t'}{(t+t')^{\alpha'}\,\sqrt{t'}}
\nonumber\\
&\stackrel{\mbox{\tiny$(t'=tr)$}}{\rule{0pt}{7pt}=}&
\|({\mathfrak A}_2\,h)'\|_{0,\alpha'}\int_0^\infty\frac{t\,\dd r}{(t+tr)^{\alpha'}\,\sqrt{tr}}
\nonumber\\
&=&
t^{-(\alpha'-\frac{1}{2})}\,
\|({\mathfrak A}_2\,h)'\|_{0,\alpha'}
\int_0^\infty\frac{\dd r}{(1+r)^{\alpha'}\,\sqrt{r}}\label{manipulation}\\
&=&
{\cal O}(t^{-(\alpha'-\frac{1}{2})}),\quad t\to\infty,
\nonumber
\end{eqnarray}
proving ${\mathfrak A}_1{\mathfrak A}_2 h\in C_{0,\alpha}$ for $0\le\alpha<2$.

We continue with the proof of item(ii) of Proposition \ref{new drift}.
First fix $h\in{\mathscr D}$ and observe that,
by (\ref{Young}) and ${\mathfrak A}_2 h\in C_{0,\alpha}$ for $0\le\alpha<3/2$,
the convolution $(4\pi|\cdot|)^{-1/2}\ast({\mathfrak A}_2 h)^a$ is well-defined and
\begin{equation}\label{principal value}
\frac{1}{\sqrt{4\pi|\cdot|}}\ast({\mathfrak A}_2 h)^a
\,=\,
\frac{1}{\sqrt{4\pi|\cdot|}}\ast
\left(\rule{0pt}{11pt}\right.
\frac{sgn(\cdot)}{\sqrt{\pi|\cdot|}}\ast(h^a)'
\left.\rule{0pt}{11pt}\right)
\,=\,h.
\end{equation}
This is easiest seen by taking the Fourier transform of
$\frac{1}{\sqrt{4\pi|\cdot|}}\ast\frac{sgn(\cdot)}{\sqrt{\pi|\cdot|}}$
which is equal to the (principal value) tempered distribution
$\tau\mapsto\frac{1}{{\bf i}\tau}$.

The next step is based on the following classical result \cite[Cor. 3.7]{dP71}
on the extension of the Stone-Weierstrass theorem for weighted topologies:
\begin{lemma}
The linear hull of $\{e^{-\nu\,\fatdot}:\nu>0\}$ is dense in the Banach space
$(C_{0,\alpha}\,,\|\cdot\|_{0,\alpha})$ for each $\alpha\ge 0$.
\end{lemma}
So, for fixed $\alpha_0\in(2,5/2)$, we can choose
$\tilde{\bf e}_n\in\LL\{e^{-\nu\,\fatdot}:\nu>0\},\,n=1,2,\dots$, such that
$\tilde{\bf e}_n\to({\mathfrak A}_2 h)'$ in $C_{0,\alpha_0}$ if $n\to\infty$.
Here, $\alpha_0<5/2$ is required for $({\mathfrak A}_2 h)'\in C_{0,\alpha_0}$
and the reason for $\alpha_0>2$ will become clear later.

Define
$${\bf e}_n(t)\,=\,-\int_t^\infty\tilde{\bf e}_n(t')\,\dd t',\quad n=1,2,\dots,$$
and observe that, for all $q\ge 1$,
\begin{eqnarray*}
&&
\int_0^\infty[\,{\bf e}_n(t)-{\mathfrak A}_2 h(t)\,]^q\,\dd t\\
&=&
\int_0^\infty[
\int_t^\infty\hspace{-5pt}
\left(\rule{0pt}{10pt}\right.
\tilde{\bf e}_n(t')-({\mathfrak A}_2 h)'(t')
\left.\rule{0pt}{10pt}\right)
\,\dd t'\,]^q\,\dd t\\
&\le&
\|\tilde{\bf e}_n-({\mathfrak A}_2 h)'\|^q_{0,\alpha_0}
\underbrace{
\int_0^\infty[
\int_t^\infty\hspace{-5pt}
(1+t')^{-\alpha_0}
\,\dd t'\,]^q\,\dd t
}_{\mbox{\tiny$<\,\infty$ because $\alpha_0>2$}}
\,\to\,0,\quad n\to\infty,
\end{eqnarray*}
hence, for all $q\ge 1$,
\begin{equation}\label{conv1}
{\bf e}_n\,\stackrel{\mbox{\tiny$L^q([0,\infty))$}}
{\rule{0pt}{7pt}\longrightarrow}\,{\mathfrak A}_2 h
\quad\mbox{and}\quad
{\bf e}'_n\,\stackrel{\mbox{\tiny$C_{0,\alpha_0}$}}
{\rule{0pt}{7pt}\longrightarrow}\,({\mathfrak A}_2 h)'
\quad\mbox{if}\quad
n\to\infty.
\end{equation}

Furthermore, for each $n$, since ${\bf e}_n\in\LL\{e^{-\nu\,\fatdot}:\nu>0\}$,
we know that $(4\pi|\cdot|)^{-1/2}\ast{\bf e}_n^a$
is a bounded continuous function in $L^2([0,\infty))$ 
by Lemma \ref{identify rho}.
But, as being shown in the next paragraph, even
\begin{equation}\label{conv2}
\frac{1}{\sqrt{4\pi|\cdot|}}\ast{\bf e}_n^a
\,\stackrel{\mbox{\tiny$L^2([0,\infty))$}}
{\rule{0pt}{7pt}\longrightarrow}\,h,\quad n\to\infty,
\end{equation}
holds true.

In fact, using (\ref{principal value}), we can write
\begin{eqnarray*}
&&\left|\rule{0pt}{11pt}\right.
(\frac{1}{\sqrt{4\pi|\cdot|}}\ast{\bf e}_n^a-h)(t)
\left.\rule{0pt}{11pt}\right|
\,=\,
\left|\rule{0pt}{11pt}\right.
\int_\NR\frac{1}{\sqrt{4\pi|t-t'|}}\,
[\,{\bf e}_n(t)-{\mathfrak A}_2 h(t)\,]^a(t')\,\dd t'
\left.\rule{0pt}{11pt}\right|\\
&=&\rule{0pt}{20pt}
\frac{1}{\sqrt{4\pi}}
\left|\rule{0pt}{11pt}\right.
\int_0^\infty(\frac{1}{\sqrt{|t-t'|}}-\frac{1}{\sqrt{t+t'}})\,
[
\int_{t'}^\infty\hspace{-5pt}
\left(\rule{0pt}{10pt}\right.
\tilde{\bf e}_n(s)-({\mathfrak A}_2 h)'(s)
\left.\rule{0pt}{10pt}\right)
\dd s\,]\,
\dd t'
\left.\rule{0pt}{11pt}\right|\\
&\le&\rule{0pt}{20pt}
\frac{1}{\sqrt{4\pi}}
\int_0^\infty(\frac{1}{\sqrt{|t-t'|}}-\frac{1}{\sqrt{t+t'}})\,
[
\int_{t'}^\infty\hspace{-5pt}
(1+s)^{-\alpha_0}\,
\dd s\,]\,
\dd t'\;
\|\tilde{\bf e}_n-({\mathfrak A}_2 h)'\|_{0,\alpha_0}\\
&\stackrel{\mbox{\tiny$(t'=tr)$}}{\rule{0pt}{7pt}=}&\rule{0pt}{20pt}
(\alpha_0-1)\;
\underbrace{
\frac{\sqrt{t}}{\sqrt{4\pi}}
\int_0^\infty(\frac{1}{\sqrt{|1-r|}}-\frac{1}{\sqrt{1+r}})\,
(1+tr)^{-\alpha_0+1}\,
\dd r
}
\;\;\|\tilde{\bf e}_n-({\mathfrak A}_2 h)'\|_{0,\alpha_0}
\end{eqnarray*}
where, if $2<\alpha_0<3$, then
the underbraced term is ${\cal O}(t^{-\alpha_0+\frac{3}{2}}),\,t\to\infty$,
by a calculation similar to how the right-hand side
of (\ref{rep l_nu}) was shown to be ${\cal O}(t^{-3/2}),\,t\to\infty$,
in the proof of Lemma \ref{identify rho}.
Therefore, when choosing $\alpha_0\in(2,5/2)$ as we do,
this large\,-$t$-behaviour of 
$[(4\pi|\cdot|)^{-1/2}\ast{\bf e}_n^a-h](t)$ can be assumed
proving (\ref{conv2}) since, for $p>2$,
\begin{eqnarray*}
&&
\left\|\rule{0pt}{11pt}\right.
\left[\rule{0pt}{10pt}\right.
\frac{1}{\sqrt{4\pi|\cdot|}}\ast{\bf e}_n^a-h
\left.\rule{0pt}{10pt}\right]\ind_{[0,1]}
\left.\rule{0pt}{11pt}\right\|_{L^2([0,\infty))}
\,=\,
\left\|\rule{0pt}{11pt}\right.
\left[\rule{0pt}{10pt}\right.
\frac{1}{\sqrt{4\pi|\cdot|}}\ast({\bf e}_n-\mathfrak{A}_2 h)^a
\left.\rule{0pt}{10pt}\right]\ind_{[0,1]}
\left.\rule{0pt}{11pt}\right\|_{L^2([0,\infty))}\\
&\le&\rule{0pt}{20pt}
\left\|\rule{0pt}{11pt}\right.
\frac{1}{\sqrt{4\pi|\cdot|}}\ast({\bf e}_n-\mathfrak{A}_2 h)^a
\left.\rule{0pt}{11pt}\right\|_{L^p(\NR)}
\,\stackrel{\mbox{\tiny(\ref{Young})}}{\rule{0pt}{10pt}\le}\,
c_p\left(\rule{0pt}{10pt}\right.
\|({\bf e}_n-\mathfrak{A}_2 h)^a\|_{L^2(\NR)}
+\|({\bf e}_n-\mathfrak{A}_2 h)^a\|_{L^1(\NR)}
\left.\rule{0pt}{10pt}\right)
\end{eqnarray*}
the right-hand side of which converges to zero by (\ref{conv1}).

Now realise
that $\sqrt{\nu}e^{-\nu\,\fatdot}={\mathfrak A}_1 e^{-\nu\,\fatdot}$ for all $\nu>0$.
Thus, applying Lemma \ref{identify rho} again, the functions 
$\varrho_{\frac{1}{\sqrt{4\pi|\cdot|}}\ast{\bf e}_n^a}$ 
exist and
$$\varrho_{\frac{1}{\sqrt{4\pi|\cdot|}}\ast{\bf e}_n^a}(U(x_0,\cdot),y)
\,\stackrel{\mbox{\tiny\rm a.s.}}{=}\,
U(y,{\mathfrak A}_1{\bf e}_n)+\partial_1 U(y,{\bf e}_n),\quad y>x_0,$$
for all $n=1,2,\dots$\,
Using the versions found in Proposition \ref{stateSpace}(v),
it follows that the process
$$\tilde{W}_z(\frac{1}{\sqrt{4\pi|\cdot|}}\ast{\bf e}_n^a)
\,\stackrel{\mbox{\tiny\rm def}}{=}\,
W_z(\frac{1}{\sqrt{4\pi|\cdot|}}\ast{\bf e}_n^a)\,-
\int_{x_0}^{x_0+z}\hspace{-5pt}[\,
U(y,{\mathfrak A}_1{\bf e}_n)+\partial_1 U(y,{\bf e}_n)
\,]\,\dd y,\;z\ge 0,$$
has all properties of a process $\{\tilde{W}_z(l);\,z\ge 0\}$
given in Proposition \ref{new wipro} 
when replacing $l$ by $\frac{1}{\sqrt{4\pi|\cdot|}}\ast{\bf e}_n^a,\,n=1,2,\dots$\,

If we can now verify that, for each $z\ge 0$, the three sequences of random variables
\begin{equation}\label{three sequences}
W_z(\frac{1}{\sqrt{4\pi|\cdot|}}\ast{\bf e}_n^a),\;
\int_{x_0}^{x_0+z}\hspace{-5pt}U(y,{\mathfrak A}_1{\bf e}_n)\,\dd y,\;
\int_{x_0}^{x_0+z}\hspace{-5pt}\partial_1 U(y,{\bf e}_n)\,\dd y,
\quad n=1,2,\dots,
\end{equation}
converge to
$$W_z(h),\;
\int_{x_0}^{x_0+z}\hspace{-5pt}U(y,{\mathfrak A}_1{\mathfrak A}_2 h)\,\dd y,\;
\int_{x_0}^{x_0+z}\hspace{-5pt}\partial_1 U(y,{\mathfrak A}_2 h)\,\dd y$$
in $L^2(\Omega)$ when $n\to\infty$, respectively,
then, for each $z\ge 0$, 
the sequence of random variables
$\tilde{W}_z(\frac{1}{\sqrt{4\pi|\cdot|}}\ast{\bf e}_n^a),\,n=1,2,\dots$,
converges in $L^2(\Omega)$
to $\tilde{W}_z(h)$ as defined in item (ii) of Proposition \ref{new drift}.
But $\{\tilde{W}_z(h);\,z\ge 0\}$ is a continuous process so, for $h\not=0$,
the process $\{\tilde{W}_z(h)/\|h\|_{L^2([0,\infty))};\,z\ge 0\}$
is a standard Wiener process with respect to the filtration $\tilde{\cal F}_z,\,z\ge 0$
by simply checking the corresponding martingale problem.
The linearity (\ref{lin new}) is an easy consequence of the properties of the summands
defining $\tilde{W}_z(h),\,z\ge 0$.

It remains to prove the convergence of the three sequences in (\ref{three sequences}).
Fix $z\ge 0$. First, the convergence
$$W_z(\frac{1}{\sqrt{4\pi|\cdot|}}\ast{\bf e}_n^a)
\,\stackrel{\mbox{\tiny$L^2(\Omega)$}}
{\rule{0pt}{7pt}\longrightarrow}\,W_z(h),\quad n\to\infty,$$
follows from (\ref{conv2}) using the definition on page \pageref{defi wipro}
of $W_z(l)$ for $l\in L^2([0,\infty))$ as a stochastic integral.

Second, applying Proposition \ref{stateSpace}(ii), we obtain that
\begin{eqnarray*}
&&\ee\,[\int_{x_0}^{x_0+z}\hspace{-5pt}
U(y,{\mathfrak A}_1{\bf e}_n-{\mathfrak A}_1{\mathfrak A}_2 h)
\,\dd y\,]^2\\
&\le&
z\int_{x_0}^{x_0+z}\hspace{-5pt}\ee\,
U(y,{\mathfrak A}_1({\bf e}_n-{\mathfrak A}_2 h))^2
\,\dd y
\,=\,
z^2\,\langle{\mathfrak A}_1({\bf e}_n-{\mathfrak A}_2 h)\,;
\frac{-\sqrt{|\cdot|}}{\sqrt{4\pi}}\ast
[{\mathfrak A}_1({\bf e}_n-{\mathfrak A}_2 h)]^a\rangle\\
&\le&
z^2\,\underbrace{
\langle(1+|\cdot|)^{-\alpha}\,;\frac{-\sqrt{|\cdot|}}{\sqrt{4\pi}}\ast
[(1+|\cdot|)^{-\alpha}]^a\rangle
}_{\mbox{\tiny$<\,\infty$ for all $\alpha>5/4$}}
\,\|{\mathfrak A}_1({\bf e}_n-{\mathfrak A}_2 h)\|^2_{0,\alpha}.
\end{eqnarray*}
So we need 
$\|{\mathfrak A}_1({\bf e}_n-{\mathfrak A}_2 h)\|_{0,\alpha}\to 0,\,n\to\infty$,
for some $\alpha>5/4$. Fix $\alpha>5/4$.
By Proposition \ref{new drift}(i), this $\alpha$ should be less than 2 to ensure that
${\mathfrak A}_1({\bf e}_n-{\mathfrak A}_2 h)\in C_{0,\alpha}$. Then
\begin{eqnarray*}
\|{\mathfrak A}_1({\bf e}_n-{\mathfrak A}_2 h)\|_{0,\alpha}
&=&
\sup_{t\ge 0}|(1+t)^\alpha\int_t^\infty
\frac{-[{\bf e}'_n(t')-({\mathfrak A}_2 h)'(t')]\,\dd t'}{\sqrt{\pi(t'-t)}}|\\
&\le&
\|{\bf e}'_n-({\mathfrak A}_2 h)'\|_{0,\alpha_0}\,
\sup_{t\ge 0}|(1+t)^\alpha\int_0^\infty
\frac{\dd t'}{(1+t'+t)^{\alpha_0}\sqrt{t'}}\,|
\end{eqnarray*}
where the last supremum is finite for $\alpha<\alpha_0-1/2$
by manipulations similar to how (\ref{manipulation}) was derived.
Recall that we can choose $\alpha\in(5/4,2)$ and $\alpha_0\in(2,5/2)$
so that the convergence of the second sequence in (\ref{three sequences}) 
follows from (\ref{conv1}).

Third, applying Proposition \ref{stateSpace}(ii) once more, we obtain that
$$\ee\,[\int_{x_0}^{x_0+z}\hspace{-5pt}
\partial_1 U(y,{\bf e}_n-{\mathfrak A}_2 h)
\,\dd y\,]^2
\,=\,
z^2\,\langle{\bf e}_n-{\mathfrak A}_2 h\,;
\frac{1}{2\sqrt{4\pi|\cdot|}}\ast({\bf e}_n-{\mathfrak A}_2 h)^a\rangle$$
where the right-hand side converges to zero 
by (\ref{principal value}),(\ref{conv1}) and (\ref{conv2})
which completes the discussion of the convergence
of the sequences in (\ref{three sequences}).

We finally show item (iii) of Proposition \ref{new drift}.
Fix $h\in{\mathscr D}$.
First, since $((4\pi|\cdot|)^{-1/2})^F=1/\sqrt{2|\cdot|}$, the equality
${\mathfrak A}_2 h\,=\,\sqrt{2}\,(-\partial_t^2)^{\frac{1}{4}}\,h^a$
follows from (\ref{principal value}) by taking Fourier transforms.
Second, observe that
$${\mathfrak A}_1 h
\,=\,
\frac{-\ind_{(-\infty,0)}}{\sqrt{\pi|\cdot|}}\ast(h^a)'
\,=\,
\frac{1}{2}
\left(\rule{0pt}{11pt}\right.
\frac{sgn(\cdot)}{\sqrt{\pi|\cdot|}}
-
\frac{1}{\sqrt{\pi|\cdot|}}
\left.\rule{0pt}{11pt}\right)
\ast(h^a)'$$
hence,
using the regularity of 
${\mathfrak A}_2 h$ as stated in Proposition \ref{new drift}(i),
the wanted equality for ${\mathfrak A}_1{\mathfrak A}_2 h$ follows from
\begin{eqnarray*}
{\mathfrak A}_1{\mathfrak A}_2 h
&=&
\frac{1}{2}{\mathfrak A}_2^2 h
\,-\,
\frac{1}{2}
\frac{1}{\sqrt{\pi|\cdot|}}\ast
\left(\rule{0pt}{11pt}\right.
({\mathfrak A}_2 h)^a
\left.\rule{0pt}{11pt}\right)'\\
&=&\rule{0pt}{20pt}
\frac{1}{2}{\mathfrak A}_2^2 h
\,-\,
\partial_t\,[\,\frac{1}{\sqrt{4\pi|\cdot|}}\ast
({\mathfrak A}_2 h)^a
\,]\\
&=&\rule{0pt}{20pt}
\frac{1}{2}{\mathfrak A}_2^2 h
\,-\,h'
\end{eqnarray*}
where, for the last line above, we again applied (\ref{principal value}).\hfill$\Box$\\

{\it Proof of} {\bf Theorem \ref{new SPDE}}.
In this proof one has to differ between the regular generalised function $U$ on $C_c^\infty({\cal Q}_+)$
given by the continuous version of the right-hand side of (\ref{greenrep})
and the version of the process $\{(U(x,\cdot),\partial U(x,\cdot));\,x\ge x_0\}$
provided by Proposition \ref{stateSpace}(v).

First fix $x_0\in\NR$ and $(f,h)\in\Co\times{\mathfrak D}$.
Then, using (\ref{equal to U}), we have that 
$U(-f'\otimes h)\stackrel{\mbox{\tiny\rm a.s.}}{=}-\int_\NR f'(x)\,U(x,h)\,\dd x$
hence, by partial integration,
the first equation of (\ref{new SDE}) yields
$U(-f'\otimes h)
\,\stackrel{\mbox{\tiny\rm a.s.}}{=}\,
\int_\NR f(x)\,\partial_1 U(x,h)\,\dd x.$
As a consequence,
$U(-f'\otimes h)=\int_\NR f(x)\,\partial_1 U(x,h)\,\dd x$
for all $(f,h)$ from a countable dense subset of $\Co\times\mathscr{D}$ almost surely.
Of course, the map
$\Co\times\mathscr{D}\to\NR:(f,h)\mapsto U(-f'\otimes h)$
is continuous. Nevertheless, by (\ref{bound on norms}), it holds that
$$\pp\left(\rule{0pt}{11pt}\right.
\int_{x_0}^x\|\partial_1 U(y,\cdot)\|_{E_2}\,\dd y\,<\,\infty\quad
\mbox{for all $x\ge x_0$}
\left.\rule{0pt}{11pt}\right)\,=\,1$$
hence, since $\mathscr{D}$ is continuously embedded in $E_2'$, the map
$\Co\times\mathscr{D}\to\NR:(f,h)\mapsto\int_\NR f(x)\,\partial_1 U(x,h)\,\dd x$
is continuous almost surely leading to
\begin{equation}\label{identity with U'}
-U(f'\otimes h)\,=\int_\NR f(x)\,\partial_1 U(x,h)\,\dd x
\end{equation}
for all $(f,h)\in\Co\times\mathscr{D}$ almost surely.

Next, when integrating both sides of the second equation of (\ref{new SDE}) by $-f'$,
we obtain that
\begin{equation}\label{draft new SPDE}
\left.\begin{array}{rcl}
U(f''\otimes h)&\stackrel{\mbox{\tiny\rm a.s.}}{=}&\displaystyle
-U(f\otimes(-\partial_t^2)^{1/2}h^a)
-\sqrt{2}\int_\NR f(x)\,\partial_1 U(x,(-\partial_t^2)^{1/4}h^a)\,\dd x\\
&&\rule{0pt}{20pt}\displaystyle\hspace{3.7cm}
+\int_\NR f'(x)\tilde{W}_{x-x_0}(h)\,\dd x
\end{array}\right\}
\end{equation}
where we have used (\ref{equal to U}), the first equation of (\ref{new SDE}) and
partial integration.
Here, by Proposition \ref{stateSpace}(iii) and Proposition \ref{new drift}(i),
the first summand on the right-hand side has the meaning of
$-\int_\NR\int_0^\infty U(x,t)\,f(x)[(-\partial_t^2)^{1/2}h^a](t)\,\dd t\dd x$.
But if we want to write $\sqrt{2}\,U(f'\otimes(-\partial_t^2)^{1/4}h^a)$ for the second summand 
then, by Remark \ref{new meaning U}, the meaning of $U$ needs to be extended.

Recall that a measurable version of the process $\{\partial U(x,\cdot);\,x\ge x_0\}$
taking values in $E_2$ was chosen at the beginning of the proof. Furthermore,
let $\Omega_0\subseteq\Omega$ be such that, on $\Omega_0$, both holds true:
$\int_{x_0}^x\|\partial_1 U(y,\cdot)\|_{E_2}\,\dd y\,<\,\infty$ for all $x\ge x_0$ 
and (\ref{identity with U'}) is satisfied for all $(f,h)\in\Co\times\mathscr{D}$.
Then, since Proposition \ref{new drift}(i) gives
$(-\partial_t^2)^{1/4}(\mathscr{D}^a)\subseteq E_2'$,
the map $\Co\to\mathscr{D}':f\mapsto\int_\NR f(x)\,\partial_1 U(x,(-\partial_t^2)^{1/4}[\cdot]^a)\,\dd x$
is continuous on $\Omega_0$ so that, for each $\omega\in\Omega_0$, the map
$$\Co\times\mathscr{D}\to\NR:(f,h)\mapsto
\int_\NR f(x)\,\partial_1 U(\omega,x,(-\partial_t^2)^{1/4}h^a)\,\dd x$$
extends to a generalised function on $C_c^\infty({\cal Q}_+^{x_0})$ by Schwartz' kernel theorem.
Also, if $\chi_N$ is a symmetric $C^\infty$-\,function on $\NR$
such that $\chi_N\equiv 1$ on $[-N,N]$ and ${\rm supp}\,\chi_N\subseteq(-N-1,N+1),\,N=1,2,\dots$,
then $\chi_N(-\partial_t^2)^{1/4}h^a\to(-\partial_t^2)^{1/4}h^a,\,N\to\infty$, in $E_2'$ thus
$$-\lim_{N\uparrow\infty}U(f'\otimes[\chi_N(-\partial_t^2)^{1/4}h^a])
\,=
\int_\NR f(x)\,\partial_1 U(x,(-\partial_t^2)^{1/4}h^a)\,\dd x$$
for all $(f,h)\in\Co\times\mathscr{D}$ on $\Omega_0$.
As a consequence, since the tensor product
$\Co\otimes\mathscr{D}$ is dense in $C_c^\infty({\cal Q}_+^{x_0})$,
the generalised function given on $\Omega_0$ by the above right-hand side
must be equal to
\begin{equation}\label{extension U}
-\lim_{N\uparrow\infty}\,(\chi_N U)(\partial_1(-\partial_2^2)^{1/4}f^a),
\quad f\in C_c^\infty({\cal Q}_+^{x_0}),\quad\mbox{on $\Omega_0$,}
\end{equation}
where $\chi_N U$ is short for $\chi_N(t)\,U(x,t),\,(x,t)\in{\cal Q}_+$.
Setting $\partial_x(-\partial_t^2)^{1/4}\,U^a$ to be the above limit on $\Omega_0$
and zero otherwise therefore defines a meaningful 
generalised function on $C_c^\infty({\cal Q}_+^{x_0})$.

Now realise that $\Omega_0$ can be chosen to be of probability one.
Repeating the above construction for $x_0>x_1>\dots$ where $x_k\to-\infty,\,k\to\infty$,
gives meaningful definitions of $\partial_x(-\partial_t^2)^{1/4}\,U^a$
based on subsets $\Omega_k$ of measure one, $k=0,1,2,\dots$\,
Of course, $\pp(\bigcap_{k\ge 0}\Omega_k)=1$ and the definition of
$\partial_x(-\partial_t^2)^{1/4}\,U^a$ based on $\bigcap_{k\ge 0}\Omega_k$
is consistent in the following sense:
if $x_{k+1}<x_k$ and ${\rm supp}\,f\subseteq{\cal Q}_+^{x_k}$ then
the limit defining $\partial_x(-\partial_t^2)^{1/4}\,U^a$ remains the same
regardless of whether it was constructed with respect to $x_{k+1}$ or $x_k$.
In this sense, $\partial_x(-\partial_t^2)^{1/4}\,U^a$ can be considered
a meaningful generalised function on $C_c^\infty({\cal Q}_+)$
which is independent of the choice of $x_0$.

Using this definition of $\partial_x(-\partial_t^2)^{1/4}\,U^a$,
equation (\ref{draft new SPDE}) becomes
$$\begin{array}{rcl}
&&\displaystyle
\partial_x^2\,U(f\otimes h)
+(-\partial_t^2)^{1/2}\,U^a(f\otimes h)
+\sqrt{2}\,\partial_x(-\partial_t^2)^{1/4}\,U^a(f\otimes h)\\
&\stackrel{\mbox{\tiny\rm a.s.}}{=}&\rule{0pt}{20pt}\displaystyle
\int_\NR f'(x)\tilde{W}_{x-x_0}(h)\,\dd x
\end{array}$$
for all $(f,h)\in\Co\times\mathscr{D}$.
Assume for a moment that we can show that one could construct 
a Brownian sheet $\tilde{B}$ on $(\Omega,{\cal F},\pp)$ such that
$$\int_\NR f'(x)\tilde{W}_{x-x_0}(h)\,\dd x
\,\stackrel{\mbox{\tiny\rm a.s.}}{=}\,
\int_\NR\int_0^\infty\tilde{B}(x,t)f'(x)h'(t)\,\dd t\dd x$$
for all $(f,h)\in\Co\times\mathscr{D}$.
Then, by continuity, the last equation can be extended to
\begin{equation}\label{final draft new SPDE}
\partial_x^2\,U(f)
+(-\partial_t^2)^{1/2}\,U^a(f)
+\sqrt{2}\,\partial_x(-\partial_t^2)^{1/4}\,U^a(f)
\,=\,
\partial_x\partial_t\tilde{B}(f)
\end{equation}
for all $f\in C_c^\infty({\cal Q}_+^{x_0})$ almost surely.
Observe that the above left-hand side does not depend on the choice of $x_0$ hence,
if the same Brownian sheet $\tilde{B}$ can be used for all $x_0$,
then (\ref{final draft new SPDE}) can easily be extended to hold
for all $f\in C_c^\infty({\cal Q}_+)$ almost surely proving the theorem.

It remains to construct the Brownian sheet $\tilde{B}$.
By Proposition \ref{new drift}, in particular using (\ref{lin new}), 
the process $\int_\NR f'(x)\tilde{W}_{x-x_0}(h)\,\dd x$
indexed by $(f,h)\in\Co\times\mathscr{D}$ is a centred Gaussian process with covariance
$$\int_\NR\int_0^\infty f_1(x)h_1(t)f_2(x)h_2(t)\,\dd t\dd x$$
which can be represented by
$$\eta(f,h)
\,\stackrel{\mbox{\tiny def}}{=}\,
\left(\rule{0pt}{11pt}\right.\!
\partial_t\,U
+(-\partial_t^2)^{1/2}\,U^a
+\sqrt{2}\,\partial_x(-\partial_t^2)^{1/4}\,U^a
\!\left.\rule{0pt}{11pt}\right)(f\otimes h)
-\int\hspace{-10pt}\int_{{\cal Q}_+}\hspace{-10pt}B(\dd x,\dd x)(f\otimes h)(x,t)$$
independently of $x_0$.
Thus, the process $\eta$ extends to a centred Gaussian process
indexed by $L^2(\NR)\times L^2([0,\infty))$ and the continuous version of the process
$\{\eta(sgn(x)\ind_{[x\wedge 0,x\vee 0]},\ind_{[0,t]});$ $\,(x,t)\in{\cal Q}_+\}$
gives the wanted Brownian sheet $\tilde{B}$.\hfill$\Box$\\

{\it Proof} of {\bf Theorem \ref{strMtheo}}.
According to the sequence of arguments laid down in the Results-Section between 
Theorem \ref{new SPDE} and Theorem \ref{strMtheo}, there is a version of the process
$\{(U(x_0+z,\cdot)$, $\partial_1 U(x_0+z,\cdot));\,z\ge 0\}$
which satisfies all conditions of Definition \ref{our mart problem}.
Hence, by Thm.4.2(b) in \cite{EK1986} and Remark \ref{improve EK theorem},
if (wp) on page \pageref{wp} holds true for the martingale problem of Definition \ref{our mart problem}
then $\{(U(x_0+z,\cdot)$, $\partial_1 U(x_0+z,\cdot));\,z\ge 0\}$ is a strong Markov process
in the sense of Definition \ref{defi strong markov}.
Moreover, it is stationary by Proposition \ref{stateSpace}(v).

The condition (wp) will be shown in two steps.
First, for an arbitrary solution $\{(u_z,v_z);$ $z\ge 0\}$ to the martingale problem
of Definition \ref{our mart problem}, we prove that
\begin{equation}\label{u-v-equation}
\left.
\begin{array}{rcl}
\langle u_z\,;h\rangle
&\stackrel{\mbox{\tiny\rm a.s.}}{=}&\displaystyle
\langle u_0\,;h\rangle+\int_0^z\langle v_y\,;h\rangle\,\dd y\\
\langle v_z\,;h\rangle
&\stackrel{\mbox{\tiny\rm a.s.}}{=}&\displaystyle\rule{0pt}{20pt}
\langle v_0\,;h\rangle-\int_0^z
\left[\rule{0pt}{11pt}\right.\!\!
\langle u_y\,;(-\partial_t^2)^{\frac{1}{2}}h^a\rangle
+\langle v_y\,;\sqrt{2}\,(-\partial_t^2)^{\frac{1}{4}}h^a\rangle
\!\!\left.\rule{0pt}{11pt}\right]
\dd y-{\mathscr W}_z(h)
\end{array}\right\}
\end{equation}
for all $(z,h)\in[0,\infty)\times\mathscr{D}$ where 
$\{{\mathscr W}_z;\,z\ge 0\}$ is a ${\mathscr D}'$\,-\,valued Wiener process
with respect to the filtration $\mathscr{F}_z$
and, second, 
we verify that the above stochastic integral equation has a pathwise unique solution
satisfying the conditions (i),(ii),(iii) of Definition \ref{our mart problem}.
This indeed shows (wp) because pathwise uniqueness of stochastic differential equations
implies weak uniqueness and weak uniqueness is sufficient 
for the uniqueness of the one-dimensional marginal distributions.

The first step is fairly standard and we only sketch the key idea.
Also, the filtration to be considered for all
martingales and Wiener processes mentioned below
is $\mathscr{F}_z,\,z\ge 0$.

Define $F_N\in{\mathfrak F}C_b^\infty(D)$ by $h\in\mathscr{D}$ and $f_N\in C_b^\infty(\NR)$
satisfying $f_N(x)=x$ for $x\in[-N,N]$ and $\sup_{x\in\NR}|f_N(x)|\le N+1,\,N=1,2,\dots$\,
Then, using both Definition \ref{our mart problem}(iv) with respect to $F_N$
and stopping times
$\inf\{z\ge 0:|\langle u_z\,;h\rangle|+|\langle v_z\,;h\rangle|\ge N\}$,
the process
$\{\langle u_z\,;h\rangle-\langle u_0\,;h\rangle-\int_0^z\langle v_y\,;h\rangle\,\dd y;\,z\ge 0\}$
can be shown to be a continuous local martingale with quadratic variation zero
which proves the first equation of (\ref{u-v-equation}).

In a similar way one shows that, for each $h\in\mathscr{D}$, the process
$\{\langle v_z\,;h\rangle-
\langle v_0\,;h\rangle+
\int_0^z
\left[\rule{0pt}{10pt}\right.
\langle u_y\,;(-\partial_t^2)^{\frac{1}{2}}h^a\rangle
+\langle v_y\,;\sqrt{2}\,(-\partial_t^2)^{\frac{1}{4}}h^a\rangle
\left.\rule{0pt}{10pt}\right]\dd y;\,z\ge 0\}$
is a continuous local martingale with quadratic variation $z\|h\|^2_{L^2([0,\infty))}$.
Here one also needs that the stochastic integral of an adapted continuous process
against a continuous local martingale always exists.

As a consequence, by the martingale characterisation of the standard Wiener process,
for each $h\in\mathscr{D}$,
there is a continuous process $\{{\mathscr W}_z(h);\,z\ge 0\}$ on $(\Omega,{\cal F},\pp)$ such that
the second equation of (\ref{u-v-equation}) is satisfied for all $z\ge 0$ almost surely
and $\{{\mathscr W}_z(h)/\|h\|_{L^2([0,\infty))};\,z\ge 0\}$ 
is a standard Wiener Process if $h\not=0$.
Of course, ${\mathscr W}_z(h)$ inherits the linearity 
$\mathscr{W}_z(a_1 h_1+a_2 h_2)\,\stackrel{\mbox{\tiny\rm a.s.}}{=}\,
a_1\mathscr{W}_z(h_1)+a_2\mathscr{W}_z(h_2)$
for each $z\ge 0,\,a_1,a_2\in\NR,\,h_1,h_2\in{\mathscr D}$
from the process $\{(u_z,v_z);$ $z\ge 0\}$ taking values in $E_1\times E_2$.
Hence, by standard theory -- see \cite{W1986} for example,
there is a version of the centred Gaussian process ${\mathscr W}_z(h)$ 
indexed by $(z,h)\in[0,\infty)\times\mathscr{D}$ such that both
the map $\mathscr{D}\to\NR:h\mapsto{\mathscr W}_z(h)$
is an element of $\mathscr{D}'$ for each $z\ge 0$ 
and the map $[0,\infty)\to\mathscr{D}':z\mapsto{\mathscr W}_z$ is continuous.
This means that $\{{\mathscr W}_z;\,z\ge 0\}$ can indeed be considered
a $\mathscr{D}'$-valued Wiener process and the first step of proving (wp) is done.

It remains to show the pathwise uniqueness of the system of equations (\ref{u-v-equation}).
So assume that two 
${\mathscr F}_z$\,-\,progressively measurable processes
$\{(u^1_z,v^1_z);$ $z\ge 0\}$ and $\{(u^2_z,v^2_z);$ $z\ge 0\}$ on $(\Omega,{\cal F},\pp)$
taking values in $E_1\times E_2$ satisfy:
\begin{itemize}
\item
$u_0^1=u_0^2$ and $v_0^1=v_0^2$;
\item
the equation (\ref{u-v-equation}) for all $(z,h)\in[0,\infty)\times\mathscr{D}$
driven by the same $\mathscr{D}'$-valued Wiener process $\{{\mathscr W}_z;\,z\ge 0\}$
with respect to the filtration $\mathscr{F}_z$ given on $(\Omega,{\cal F},\pp)$;
\item
the conditions (i),(ii),(iii) of Definition \ref{our mart problem}.
\end{itemize}
Set $u\stackrel{\mbox{\tiny def}}{=}u^1-u^2$ and $v\stackrel{\mbox{\tiny def}}{=}v^1-v^2$
and realise that
$$\begin{array}{rcl}
\langle u_z\,;h\rangle
&\stackrel{\mbox{\tiny\rm a.s.}}{=}&\displaystyle
\int_0^z\langle v_y\,;h\rangle\,\dd y\\
\langle v_z\,;h\rangle
&\stackrel{\mbox{\tiny\rm a.s.}}{=}&\displaystyle\rule{0pt}{20pt}
-\int_0^z
\left[\rule{0pt}{11pt}\right.\!\!
\langle u_y\,;(-\partial_t^2)^{\frac{1}{2}}h^a\rangle
+\langle v_y\,;\sqrt{2}\,(-\partial_t^2)^{\frac{1}{4}}h^a\rangle
\!\!\left.\rule{0pt}{11pt}\right]
\dd y
\end{array}$$
for all $(z,h)\in[0,\infty)\times\mathscr{D}$. 
We want to show that $u\equiv 0$ almost surely.

First,
by the same principles applied in the proof of Theorem \ref{new SPDE},
one can justify that
\begin{equation}\label{justify}
u\left(\rule{0pt}{11pt}\right.
\partial_z^2 f_{|_{{\cal Q}_+^0}}
+(-\partial_t^2)^{\frac{1}{2}}[f_{|_{{\cal Q}_+^0}}]^a
-\sqrt{2}\,\partial_z(-\partial_t^2)^{\frac{1}{4}}[f_{|_{{\cal Q}_+^0}}]^a
\left.\rule{0pt}{11pt}\right)
\,=\,0
\end{equation}
for all $f\in C_c^\infty({\cal Q}_+)$ almost surely
where, in this context, $u$ stands for the regular generalised function
given by $u_z(t),\,(z,t)\in{\cal Q}_+^0$,
and $f_{|_{{\cal Q}_+^0}}$ denotes the restriction of $f$ to ${\cal Q}_+^0$.
Notice that, because $f_{|_{{\cal Q}_+^0}}$ does not have compact support 
in ${\cal Q}_+^0$ for general $f\in C_c^\infty({\cal Q}_+)$,
one also has to approximate $f_{|_{{\cal Q}_+^0}}$ by functions from $C_c^\infty({\cal Q}_+^0)$
when showing (\ref{justify}).

Second,
since the map $(z,t)\mapsto u_z(t)$ is continuous on the closure of ${\cal Q}_+^0$
and zero on the boundary of ${\cal Q}_+^0$,
$$\underline{u}(z,t)\,\stackrel{\mbox{\tiny def}}{=}
\left\{\begin{array}{rcl}
0&:&z<0,\,t\in\NR\\
u(z,t)&:&z\ge 0,\,t\ge 0\\
-u(z,-t)&:&z\ge 0,\,t<0
\end{array}\right.$$
defines a continuous function on $\NR^2$. Furthermore, for $f\in C_c^\infty({\cal Q}_+)$,
we have that
$$\int_\NR\int_\NR\underline{u}(z,t)\,f^a(z,t)\,\dd t\dd z
\,=\,
\int_0^\infty 2\int_0^\infty{u}(z,t)\,f_{|_{{\cal Q}_+^0}}(z,t)\,\dd t\dd z$$
because $\underline{u}$ is also antisymmetric in $t$.
Hence (\ref{justify}) implies 
$$\underline{u}\left(\rule{0pt}{11pt}\right.
\partial_z^2 f^a
+(-\partial_t^2)^{\frac{1}{2}}f^a
-\sqrt{2}\,\partial_z(-\partial_t^2)^{\frac{1}{4}}f^a
\left.\rule{0pt}{11pt}\right)
\,=\,0$$
for all $f\in C_c^\infty({\cal Q}_+)$ almost surely.

Note that all arguments leading to the above equality remain valid
if the corresponding functions $f$ are complex-valued and all our test function spaces
are supposed to be complex-valued for the rest of this proof. 

Therefore, because $\underline{u}$ is a continuous linear form 
on the space $\mathscr{S}(\NR^2)$ of rapidly decreasing functions
and $\{f^a:f\in C_c^\infty({\cal Q}_+)\}$ is dense in $\mathscr{S}^a(\NR^2)$,
we even obtain that
$$\underline{u}\left(\rule{0pt}{11pt}\right.
\partial_z^2 f
+(-\partial_t^2)^{\frac{1}{2}}f
-\sqrt{2}\,\partial_z(-\partial_t^2)^{\frac{1}{4}}f
\left.\rule{0pt}{11pt}\right)
\,=\,0$$
for all $f\in\mathscr{S}^a(\NR^2)$ almost surely.
But each $f\in\mathscr{S}(\NR^2)$ can be split in a unique way
into a sum of two functions in $\mathscr{S}(\NR^2)$, 
one being symmetric and the other being antisymmetric in the second argument,
and $\underline{u}$ maps the symmetric one to zero.
Moreover, when a fractional Laplacian with respect to the second argument is applied to 
a function in $\mathscr{S}(\NR^2)$ which is symmetric in the second argument,
then the outcome is still symmetric in the second argument.
So, the equality
$$\underline{u}\left(\rule{0pt}{11pt}\right.
\partial_z^2 f
+(-\partial_t^2)^{\frac{1}{2}}f
-\sqrt{2}\,\partial_z(-\partial_t^2)^{\frac{1}{4}}f
\left.\rule{0pt}{11pt}\right)
\,=\,0$$
must eventually hold for all $f\in\mathscr{S}(\NR^2)$ almost surely
and we can perform the calculation
\begin{eqnarray*}
0&=&
\underline{u}\left(\rule{0pt}{11pt}\right.
\partial_z^2 f
+(-\partial_t^2)^{\frac{1}{2}}f
-\sqrt{2}\,\partial_z(-\partial_t^2)^{\frac{1}{4}}f
\left.\rule{0pt}{11pt}\right)\\
&=&\rule{0pt}{20pt}
\underline{u}^F\left(\rule{0pt}{11pt}\right.
(\partial_z^2 f)^{F^{-1}}
+\,((-\partial_t^2)^{\frac{1}{2}}f)^{F^{-1}}
-\,\sqrt{2}\,(\partial_z(-\partial_t^2)^{\frac{1}{4}}f)^{F^{-1}}
\left.\rule{0pt}{11pt}\right)\\
&=&\rule{0pt}{20pt}
\underline{u}^F\left(\rule{0pt}{11pt}\right.
-\zeta^2 f^{F^{-1}}
+\,|\tau|f^{F^{-1}}
+\,\sqrt{2}\,{\bf i}\zeta\sqrt{|\tau|}f^{F^{-1}}
\left.\rule{0pt}{11pt}\right)
\end{eqnarray*}
leading to the well-defined equality
$$\left(\rule{0pt}{10pt}\right.
-\zeta^2 +|\tau|+\sqrt{2}\,{\bf i}\zeta\sqrt{|\tau|}
\left.\rule{0pt}{10pt}\right)
\underline{u}^F
\,\equiv\,0\quad\mbox{a.s.}$$

Observe that the only solution to the equation
$\zeta^2 -|\tau|-\sqrt{2}\,{\bf i}\zeta\sqrt{|\tau|}=0$
is $\zeta=\tau=0$ hence the tempered distribution $\underline{u}^F$
must almost surely coincide with a series expansion of type
$$\sum_{\gamma_1=0}^\infty\sum_{\gamma_2=0}^\infty
c_{\gamma_1,\gamma_2}\,\partial_1^{\gamma_1}\partial_2^{\gamma_2}\,\delta_{(0,0)}$$
where $\delta_{(0,0)}$ denotes Dirac's delta-function with respect to $(0,0)\in\NR^2$
and $c_{\gamma_1,\gamma_2}$ are complex-valued coefficients.
Taking the inverse Fourier transform gives
$$\underline{u}(z,t)
\,=\,
\sum_{\gamma_1=0}^\infty\sum_{\gamma_2=0}^\infty
\frac{(-{\bf i})^{\gamma_1+\gamma_2}}{(2\pi)^2}\,
c_{\gamma_1,\gamma_2}\,z^{\gamma_1}t^{\gamma_2}
\,=\,
\sum_{\gamma_1=0}^\infty\sum_{\gamma_2=0}^\infty
\underline{c}_{\gamma_1,\gamma_2}\,z^{\gamma_1}t^{\gamma_2}$$
for all $(z,t)\in\NR^2$ almost surely
with real\,-valued coefficients $\underline{c}_{\gamma_1,\gamma_2}$
because $\underline{u}$ is real\,-valued.
But this proves $\underline{u}\equiv 0$ almost surely, 
hence the pathwise uniqueness of (\ref{u-v-equation}),
because, on the one hand, the power series 
$\sum_{\gamma_1=0}^\infty\sum_{\gamma_2=0}^\infty
\underline{c}_{\gamma_1,\gamma_2}\,z^{\gamma_1}t^{\gamma_2}$
cannot depend on $z$ as $\underline{u}(z,t)=0$ whenever $z<0$
and, on the other hand, if 
$\underline{u}(z,t)=\sum_{\gamma_2=0}^\infty\underline{c}_{\gamma_2}\,t^{\gamma_2}$
then all coefficients $\underline{c}_{\gamma_2}$ must vanish
since $\underline{u}(0,t)=0$ for all $t\in\NR$.\hfill$\Box$\\

{\it Proof} of {\bf Corollary \ref{corMtheo}}.
One only has to justify why, for fixed $x\in\NR$, the $\sigma$-algebra
$germ\left(\rule{0pt}{10pt}\right.
\{x\}\times(0,\infty)
\left.\rule{0pt}{10pt}\right)$
strictly contains
$\sigma(U(x,\cdot),\partial_1 U(x,\cdot))$.
But, the former includes information on all 
(possibly generalised)
derivatives $\partial_1^m U(x,\cdot),\,m=0,1,2,\dots$, 
while the latter does not.\hfill$\Box$

\end{document}